\documentclass[reqno,12pt,a4paper]{amsart}

\usepackage{mathrsfs}
\usepackage{amssymb}
\usepackage{amsfonts}
\usepackage{latexsym}
\usepackage{amsthm}

\usepackage{graphicx}
\def\lb{\label}

\newcommand{\er}[1]{\textrm{(\ref{#1})}}

\begin{document}

%%%%%%%%%% Some definitions %%%%%%%%%%

%%%%%%%% Equations, theorems %%%%%%%%%
\renewcommand{\theequation}{\arabic{section}.\arabic{equation}}
\theoremstyle{plain}
\newtheorem{theorem}{\bf Theorem}[section]
\newtheorem{lemma}[theorem]{\bf Lemma}
\newtheorem{corollary}[theorem]{\bf Corollary}
\newtheorem{proposition}[theorem]{\bf Proposition}
\newtheorem{definition}[theorem]{\bf Definition}
\newtheorem{remark}[theorem]{\it Remark}
%\theoremstyle{remark}
%\newtheorem{remark}[theorem]{\bf Remark}

%%%%% Alphabet %%%%%
\def\a{\alpha}  \def\cA{{\mathcal A}}     \def\bA{{\bf A}}  \def\mA{{\mathscr A}}
\def\b{\beta}   \def\cB{{\mathcal B}}     \def\bB{{\bf B}}  \def\mB{{\mathscr B}}
\def\g{\gamma}  \def\cC{{\mathcal C}}     \def\bC{{\bf C}}  \def\mC{{\mathscr C}}
\def\G{\Gamma}  \def\cD{{\mathcal D}}     \def\bD{{\bf D}}  \def\mD{{\mathscr D}}
\def\d{\delta}  \def\cE{{\mathcal E}}     \def\bE{{\bf E}}  \def\mE{{\mathscr E}}
\def\D{\Delta}  \def\cF{{\mathcal F}}     \def\bF{{\bf F}}  \def\mF{{\mathscr F}}
\def\c{\chi}    \def\cG{{\mathcal G}}     \def\bG{{\bf G}}  \def\mG{{\mathscr G}}
\def\z{\zeta}   \def\cH{{\mathcal H}}     \def\bH{{\bf H}}  \def\mH{{\mathscr H}}
\def\e{\eta}    \def\cI{{\mathcal I}}     \def\bI{{\bf I}}  \def\mI{{\mathscr I}}
\def\p{\psi}    \def\cJ{{\mathcal J}}     \def\bJ{{\bf J}}  \def\mJ{{\mathscr J}}
\def\vT{\Theta} \def\cK{{\mathcal K}}     \def\bK{{\bf K}}  \def\mK{{\mathscr K}}
\def\k{\kappa}  \def\cL{{\mathcal L}}     \def\bL{{\bf L}}  \def\mL{{\mathscr L}}
\def\l{\lambda} \def\cM{{\mathcal M}}     \def\bM{{\bf M}}  \def\mM{{\mathscr M}}
\def\L{\Lambda} \def\cN{{\mathcal N}}     \def\bN{{\bf N}}  \def\mN{{\mathscr N}}
\def\m{\mu}     \def\cO{{\mathcal O}}     \def\bO{{\bf O}}  \def\mO{{\mathscr O}}
\def\n{\nu}     \def\cP{{\mathcal P}}     \def\bP{{\bf P}}  \def\mP{{\mathscr P}}
\def\r{\rho}    \def\cQ{{\mathcal Q}}     \def\bQ{{\bf Q}}  \def\mQ{{\mathscr Q}}
\def\s{\sigma}  \def\cR{{\mathcal R}}     \def\bR{{\bf R}}  \def\mR{{\mathscr R}}
\def\S{\Sigma}  \def\cS{{\mathcal S}}     \def\bS{{\bf S}}  \def\mS{{\mathscr S}}
\def\t{\tau}    \def\cT{{\mathcal T}}     \def\bT{{\bf T}}  \def\mT{{\mathscr T}}
\def\f{\phi}    \def\cU{{\mathcal U}}     \def\bU{{\bf U}}  \def\mU{{\mathscr U}}
\def\F{\Phi}    \def\cV{{\mathcal V}}     \def\bV{{\bf V}}  \def\mV{{\mathscr V}}
\def\P{\Psi}    \def\cW{{\mathcal W}}     \def\bW{{\bf W}}  \def\mW{{\mathscr W}}
\def\o{\omega}  \def\cX{{\mathcal X}}     \def\bX{{\bf X}}  \def\mX{{\mathscr X}}
\def\x{\xi}     \def\cY{{\mathcal Y}}     \def\bY{{\bf Y}}  \def\mY{{\mathscr Y}}
\def\X{\Xi}     \def\cZ{{\mathcal Z}}     \def\bZ{{\bf Z}}  \def\mZ{{\mathscr Z}}
\def\O{\Omega}

\newcommand{\gA}{\mathfrak{A}}
\newcommand{\gB}{\mathfrak{B}}
\newcommand{\gC}{\mathfrak{C}}
\newcommand{\gD}{\mathfrak{D}}
\newcommand{\gE}{\mathfrak{E}}
\newcommand{\gF}{\mathfrak{F}}
\newcommand{\gG}{\mathfrak{G}}
\newcommand{\gH}{\mathfrak{H}}
\newcommand{\gI}{\mathfrak{I}}
\newcommand{\gJ}{\mathfrak{J}}
\newcommand{\gK}{\mathfrak{K}}
\newcommand{\gL}{\mathfrak{L}}
\newcommand{\gM}{\mathfrak{M}}
\newcommand{\gN}{\mathfrak{N}}
\newcommand{\gO}{\mathfrak{O}}
\newcommand{\gP}{\mathfrak{P}}
\newcommand{\gQ}{\mathfrak{Q}}
\newcommand{\gR}{\mathfrak{R}}
\newcommand{\gS}{\mathfrak{S}}
\newcommand{\gT}{\mathfrak{T}}
\newcommand{\gU}{\mathfrak{U}}
\newcommand{\gV}{\mathfrak{V}}
\newcommand{\gW}{\mathfrak{W}}
\newcommand{\gX}{\mathfrak{X}}
\newcommand{\gY}{\mathfrak{Y}}
\newcommand{\gZ}{\mathfrak{Z}}

\def\ve{\varepsilon}   \def\vt{\vartheta}    \def\vp{\varphi}    \def\vk{\varkappa}

\def\Z{{\mathbb Z}}    \def\R{{\mathbb R}}   \def\C{{\mathbb C}}    \def\K{{\mathbb K}}
\def\T{{\mathbb T}}    \def\N{{\mathbb N}}   \def\dD{{\mathbb D}}

%%%%% Arrows %%%%%

\def\la{\leftarrow}              \def\ra{\rightarrow}            \def\Ra{\Rightarrow}
\def\ua{\uparrow}                \def\da{\downarrow}
\def\lra{\leftrightarrow}        \def\Lra{\Leftrightarrow}

%%%%% Typography %%%%%

\def\lt{\biggl}                  \def\rt{\biggr}
\def\ol{\overline}               \def\wt{\widetilde}
\def\no{\noindent}

%%%%% Math signs %%%%%

\let\ge\geqslant                 \let\le\leqslant
\def\lan{\langle}                \def\ran{\rangle}
\def\/{\over}                    \def\iy{\infty}
\def\sm{\setminus}               \def\es{\emptyset}
\def\ss{\subset}                 \def\ts{\times}
\def\pa{\partial}                \def\os{\oplus}
\def\om{\ominus}                 \def\ev{\equiv}
\def\iint{\int\!\!\!\int}        \def\iintt{\mathop{\int\!\!\int\!\!\dots\!\!\int}\limits}
\def\el2{\ell^{\,2}}             \def\1{1\!\!1}
\def\sh{\sharp}
\def\wh{\widehat}
\def\bs{\backslash}
%%%%% Math operations %%%%%

\def\all{\mathop{\mathrm{all}}\nolimits}
\def\Area{\mathop{\mathrm{Area}}\nolimits}
\def\arg{\mathop{\mathrm{arg}}\nolimits}
\def\const{\mathop{\mathrm{const}}\nolimits}
\def\det{\mathop{\mathrm{det}}\nolimits}
\def\diag{\mathop{\mathrm{diag}}\nolimits}
\def\diam{\mathop{\mathrm{diam}}\nolimits}
\def\dim{\mathop{\mathrm{dim}}\nolimits}
\def\dist{\mathop{\mathrm{dist}}\nolimits}
\def\Im{\mathop{\mathrm{Im}}\nolimits}
\def\Iso{\mathop{\mathrm{Iso}}\nolimits}
\def\Ker{\mathop{\mathrm{Ker}}\nolimits}
\def\Lip{\mathop{\mathrm{Lip}}\nolimits}
\def\rank{\mathop{\mathrm{rank}}\limits}
\def\Ran{\mathop{\mathrm{Ran}}\nolimits}
\def\Re{\mathop{\mathrm{Re}}\nolimits}
\def\Res{\mathop{\mathrm{Res}}\nolimits}
\def\res{\mathop{\mathrm{res}}\limits}
\def\sign{\mathop{\mathrm{sign}}\nolimits}
\def\span{\mathop{\mathrm{span}}\nolimits}
\def\supp{\mathop{\mathrm{supp}}\nolimits}
\def\Tr{\mathop{\mathrm{Tr}}\nolimits}
\def\BBox{\hspace{1mm}\vrule height6pt width5.5pt depth0pt \hspace{6pt}}
\def\where{\mathop{\mathrm{where}}\nolimits}
\def\as{\mathop{\mathrm{as}}\nolimits}

%%%%%%%%%%%%% specialities %%%%%%%%%%%%%%

\newcommand\nh[2]{\widehat{#1}\vphantom{#1}^{(#2)}}
%{{\mathop{#1}\limits^\wedge}\vphantom{#1}^{(#2)}}
\def\dia{\diamond}

\def\Oplus{\bigoplus\nolimits}

%%%%%%%%%%% End of definitions %%%%%%%%%%

%%%%% OLD OLD OLD

\def\qqq{\qquad}
\def\qq{\quad}
\let\ge\geqslant
\let\le\leqslant
\let\geq\geqslant
\let\leq\leqslant
\newcommand{\ca}{\begin{cases}}
\newcommand{\ac}{\end{cases}}
\newcommand{\ma}{\begin{pmatrix}}
\newcommand{\am}{\end{pmatrix}}
\renewcommand{\[}{\begin{equation}}
\renewcommand{\]}{\end{equation}}
\def\eq{\begin{equation}}
\def\qe{\end{equation}}
\def\[{\begin{equation}}
\def\bu{\bullet}

\def\Wr{\mathop{\rm Wr}\nolimits}
\def\BBox{\hspace{1mm}\vrule height6pt width5.5pt depth0pt \hspace{6pt}}

\def\Diag{\mathop{\rm Diag}\nolimits}

\title[{Zigzag nanoribbons in external fields}]
        {Zigzag nanoribbons in external electric and magnetic fields}
\date{\today}

\date{\today}
\author[Evgeny L. Korotyaev]{Evgeny L. Korotyaev}
\address{School of Math., Cardiff University.
Senghennydd Road, CF24 4AG Cardiff, Wales, UK.  \ Current address:
Saint-Petersburg State University of Technology and Design, Bolshaya
Morskaya 18, Russia, e-mail: korotyaev@gmail.com, {\rm Partially
supported by EPSRC grant EP/D054621.}}

\author[Anton A. Kutsenko]{Anton A. Kutsenko}
\address{Laboratoire de M\'ecanique Physique, UMR CNRS 5469,
Universit\'e Bordeaux 1, Talence 33405, France,  \qqq email \
kucenkoa@rambler.ru }


\maketitle

\begin{abstract}
\no We consider the Schr\"odinger operators on zigzag nanoribbons
 (quasi-1D tight-binding models) in external magnetic fields and
an electric potential  $V$. The magnetic field is perpendicular to
the plane of the ribbon and the electric field is perpendicular to
the axis of the nanoribbon and the magnetic field. If the magnetic
and electric fields are absent, then the spectrum of the
Schr\"odinger (Laplace) operator consists of two non-flat bands and
one flat band (an eigenvalue with infinite multiplicity) between
them. If we switch on the magnetic field, then the spectrum
of the magnetic Schr\"odinger operator consists of some non-flat
bands and one flat band between them. Thus the magnetic field
changes the continuous spectrum but does not the flat band.  If we
switch on a weak electric potential $V\to 0$, then there are two
cases: (1) the flat band splits into the small spectral band. We
determine the asymptotics of the spectral bands for small fields.
(2) the unperturbed flat band remains the flat band. We describe all
potentials when the unperturbed flat band remains the flat band and
when one splits into the small band of the continuous spectrum.
Moreover, we solve inverse spectral problems for small potentials.

\end{abstract}

%  \vskip 0.25cm

\section {Introduction}
\setcounter{equation}{0}

\begin{figure}[t]
\label{f001}
         \centering\includegraphics[clip]{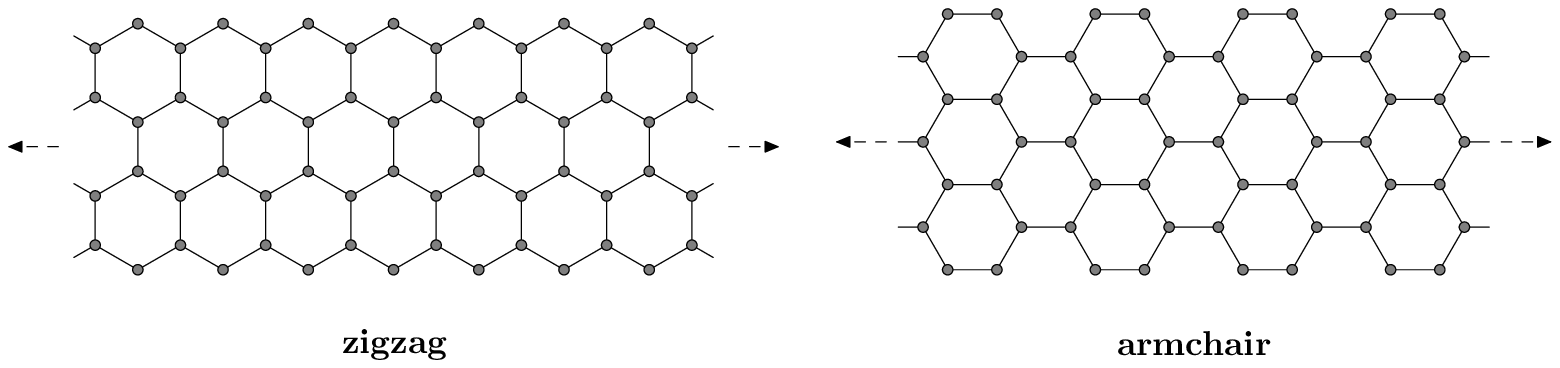}
         \begin{center}
         {\small Fig 1. Two types of horizontal nanoribbons  }
         \end{center}
         \end{figure}

There are a lot of papers about the electronic structure of carbon
materials such as carbon nanoribbons, nanotubes  and fullerenes
during the past two decades because of fundamental scientific
interest in nanomaterials and because of their versatile electronic
properties that are expected to be important for future
nanoelectronics \cite{DDE}, \cite{Ha}, \cite{SDD}. Among the carbon
nanostructures, a simple variation of graphene, ribbons has been
studied extensively. The graphene nanoribbons with varying widths
can be realized either by cutting \cite{Hi} mechanically exfoliated
graphenes \cite{No1}, or by patterning epitaxially grown graphenes
\cite{ZTSK}.

In our paper we concentrate on carbon nanoribbons, the recently
discovered two-dimensional carbon crystal \cite{No1}, \cite{No2}.
The  graphene forms a hexagonal graph (the honeycomb lattice)
embedded in $\R^2$ with all edges of constant length, see Fig.1-3.
In physics one considers only two type of ribbons: zigzag  and
armchair. When a single graphite layer is terminated by zigzag edges
on both sides, then this layer is a zigzag graphene nanoribbon
(ZGNR), see Fig.1. When a single graphite layer is terminated by
armchair edges on both sides, then this layer is a armchair graphene
nanoribbon, see Fig.1. Following conventional notation, a ZGNR is
specified by the number $N$ of  zigzag chains along the ribbon
forming the width and is referred as a N-ZGNR. For example, in Fig.1
we have the 4-ZGNR.

The main interest in physics is to study electronic structure of
ribbons, see \cite{DDE}, \cite{Ha}, \cite{SDD} and references
therein.  Moreover, there are numerous physical papers devoted to
graphene nanotubes and nanoribons in an external field, see
\cite{KLYH}, \cite{N}, \cite{No3}, \cite{SCL}, \cite{SDD} and
references therein.

There are some mathematical papers about graphene-media, see
\cite{KL}, \cite{KL1}, \cite{K}, \cite{KuP} and \cite{Pa}, devoted
to continuous models, see also \cite{KS}, \cite{RR} about other
models. But in the physical literature the most commonly used model
is the tight-binding model, corresponding to discrete Hamiltonians.
The zigzag and armchair nanotubes in external electric and magnetic
fields (discrete models) were considered in \cite{KKu4}. The zigzag
nanoribbons  (discrete models) in external electric field were
considered in \cite{KKu3}.

We consider now the N-ZGNR, $N\ge 1$ in external fields, see Fig. 2.
In fact we have a 2D problem. Firstly, there is a 2D nanoribbon $\ss
\R^2$. The ZGNRs are assumed to be infinite along the x direction
(horizontal). In general, the boundary of the ZGNR can be very
complicated. We consider only the ZGNRs with the simple boundary,
see Fig. 2. Secondly, an external 2-dim transverse electric field
$\cE_{ext}\in \R^2$ is applied across the ZGNR and along the y
direction (vertical), see Fig. 2 with the 3-ZGNR in an external
electric potential. Moreover, the homogeneous magnetic field is
perpendicular to the ribbon and along the z direction.

We will show that the spectrum of the Laplacian (an unperturbed
operator) on  the ZNRB consists of two spectral bands separated by a
gap and exactly one eigenvalue with infinite multiplicity (a flat
band) in this gap. If we switch on an external magnetic field, then
we will show that the unperturbed spectral bands change, but the
flat band (unperturbed eigenvalue) will be still a flat band. If
additionally we switch on an small external electric field, then we
will show that the spectral bands slightly change and the flat band
will splits into either a new spectral band or will be still a flat
band. Our main goal is to describe this transformation. For
applications of our model see ref. in \cite{Ha}, \cite{SDD}. The
second part of our paper is devoted to the inverse problems.

Denote a N-ZGNR by $\G=\G_N$. We consider the ribbon $\G\ss \R^2$ as
graph. This graph is just a discrete set of vertices $\bV=\{{\bold
v}_\vk: \vk=(n,k)\in \Z\ts\N_p\}$, where
 $p=2N+1, \N_N=\{1,..,N\}$ and a set of undirected edges $\bE$
 such that $[{\bold v},{\bold v}']\in \bE$ if the vertices
 ${\bold v}$ and ${\bold v}'$ are connected by a edge.
 Each vertex inside the ribbon   is connected
with some other 3  vertices, the vertex on the boundary is connected
with 1 or 2  other vertices. If the vertices ${\bold v}$ and ${\bold
v}'$ are connected by an edge, we  denote this by ${\bold v} \sim
{\bold v}'$. Introduce  the discrete Hilbert space $\ell^2(\G)$
consisting of functions $f=(f_\vk)_{\vk\in \Z\ts \N_p}$ on the set
of vertices  $\bV$ equipped with the norm $\|f\|_{\ell^2(\G)}^2=\sum
|f_{\vk}|^2$. We define the magnetic Hamiltonian $\D_b$ on the
nanoribbon $\G$  in an external magnetic fields $\mB=B(0,0,1)\in
\R^3$ by
$$
(\D_b f)_{\vk}=\sum _{{\bold v}_\vk\sim {\bold v}_{\vk'}}
e^{ia_{{\vk},{\vk'}}}f_{\vk'},\qqq \vk=(n,k)\in \Z\ts\N_p,\qqq
p=2N+1,\qq b={B\sqrt3\/2},
$$
where   $f=(f_{\vk})_{\vk\in \Z\ts \N_p}\in\ell^2(\G)$ satisfies the
Dirichlet boundary conditions
\[
 \lb{dbc}
 f_{n,0}=f_{n,p}=0, \qq  \ n\in \Z.
\]
The  factor $e^{ia_{{\vk},{\vk'}}}$ is associated with the magnetic
field $\mB$ and if $\mB=0$, then $e^{ia_{{\vk},{\vk'}}}=1$.

We define the discrete Hamiltonian $H_b=\D_b+V$ on the nanoribbon
$\G$ (a tight-binding model of single-wall nanoribbons, see
\cite{SDD}, \cite{N}) in an external electric potential $V$ and the
uniform magnetic field $\mB=B(0,0,1)\in \R^3$. The magnetic field is
perpendicular to the plane of the ribbon and the electric field is
perpendicular to the axis of the nanoribbon and the magnetic field.

 Our model
nanoribbon $\G$ is a graph, which is a set of vertices ${\bold
v}_{n,k}$ and edges $E_{n,k,j}$  given by
\[
 \ca{\bold v}_{n,2k+1}=(\sqrt3(2n+k),3k), & k\in\N^0_N\\
 {\bold v}_{n,2k}=(\sqrt3(2n+k),3k-2), & k\in\N_{N}\ac,\ \ \ca
 E_{n,k,1}=[{\bold v}_{n,2k},{\bold v}_{n,2k+1}]\\
 E_{n,k,2}=[{\bold v}_{n,2k},{\bold v}_{n,2k-1}]\\
 E_{n,k,3}=[{\bold v}_{n,2k},{\bold v}_{n+1,2k-1}]
 \ac,\ \ n\in\Z,
\]
where  $\N_k=\{1,..,k\}\ss\N$ and $\N^0_k=\N_k\cup\{0\}$, see Fig. 2
and 3 for the case $N=3$.
\begin{figure}[t]
\label{f002}
         \centering\includegraphics[clip]{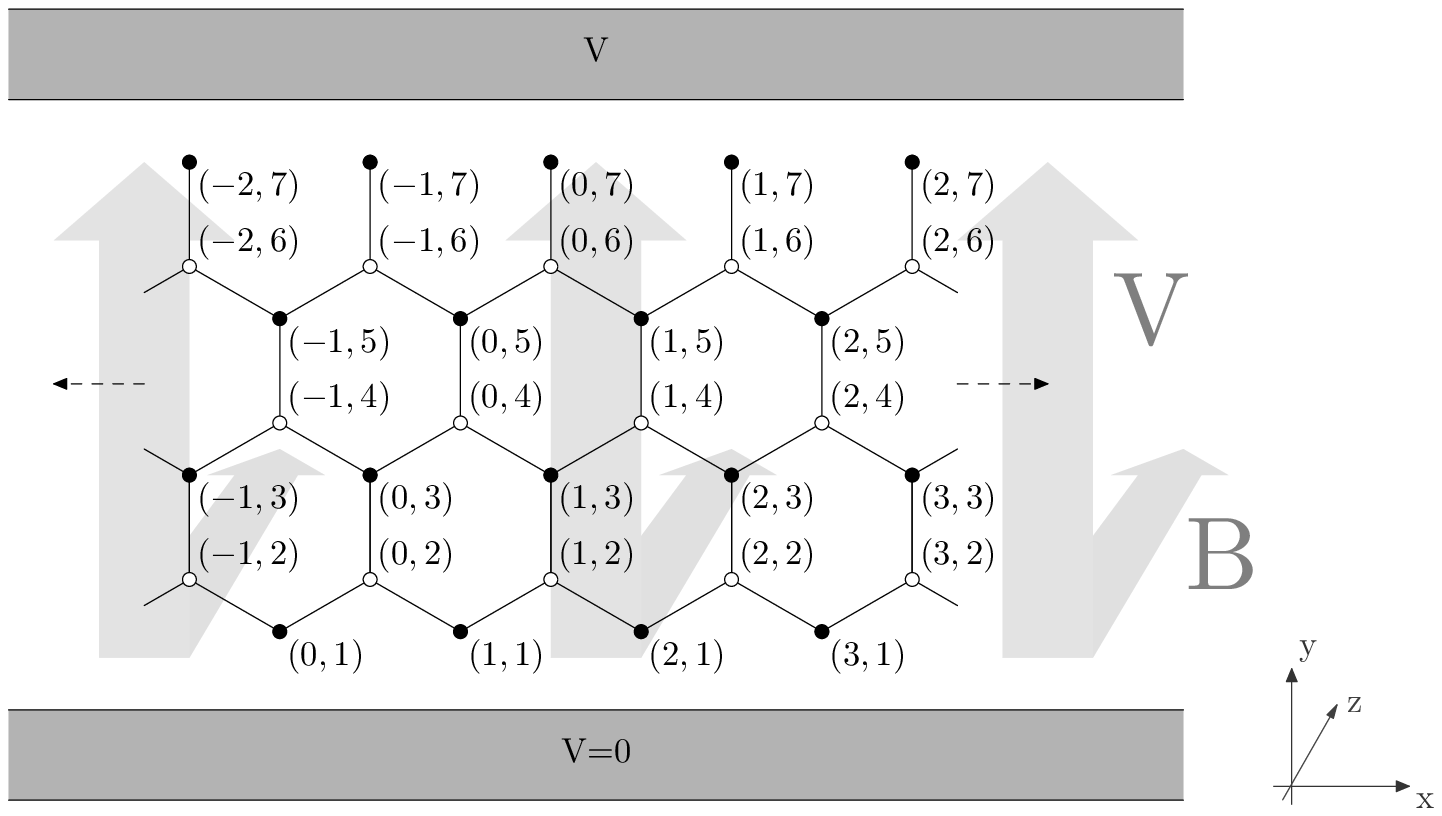}
         \begin{center}
         {\small Fig 2. A  horizontal zigzag nanoribbon at $N=3, k\in \{1,2,..,7\}$ in the constant vertical electric field. }
         \end{center}
         \end{figure}

Our Schr\"odinger operator has the form $H_b=\D_b+V,
b={B\sqrt3\/2}\in \T=\R/(2\pi \Z)$, where the magnetic operator
$\D_b$ (see more in Section 4) and the electric potential
 $V$ are given by
\begin{multline}
 \lb{H}
 (\D_bf)_{n,2k+1}=
 e^{-ib(2n+k)}f_{n,2k}+e^{ib(n+2k)}f_{n-1,2k+2}+e^{ib(n-k)}f_{n,2k+2},\qq
 k\in \N^0_{N}, \\
 f_{n,0}=f_{n,p+1}=0, \qqq p=2N+1,
 \\
 (\D_bf)_{n,2k}=e^{-ib(n-k+1)}f_{n,2k-1}+e^{-ib(n+2k-1)}f_{n+1,2k-1}
 +e^{ib(2n+k)}f_{n,2k+1},
 \ \ k\in\N_N,
\end{multline}
\[
\lb{V} (Vf)_\vk=v_k f_\vk,\qqq \vk=(n,k)\in \Z\ts\N_{p},\qq
\]
where  $n\in \Z$ and $f=(f_\vk)_{\vk\in \Z\ts \N_p}\in\ell^2(\G)$
and $v=(v_k)_1^p\in\R^p$. In fact we consider the Schr\"odinger
operator $H_b$ on the ribbon $\G$ (on the set $\Z\ts\N_{p}$) with
the Dirichlet boundary conditions $f_{n,0}=f_{n,p+1}=0$. Note that
\er{H} gives $H_{b+2\pi}=H_{b}$ for all $b\in \R$. Our electric
potential $V$ is given by $v\in\R^p$, since the electric field is
perpendicular to the axis of the nanoribbon. We formulate our
preliminary result.

\begin{theorem}
\label{T1} i)  The operator $H_b=\D_b+V$ is unitarily equivalent to
the operator $\int_{[0,2\pi)}^{\os}J_t{dt\/2\pi }$, where $J_t\ev
J_t(b,v)$  is a  Jacobi operator, acting on $\C^{p}$ and given by
\[
 \lb{Jk}
 (J_t y)_n=a_{n-1}y_{n-1}+a_{n}y_{n+1}+v_ny_n,\qq
a_{n}=\ca  1 & {\rm even}\  n \\
2|\cos(\frac t2-\frac{3n-2}2 b)|& {\rm odd}  \ n\ac,
\]
%end{multline}
where $y=(y_n)_1^{p}\in\C^{p},\qq y_0=0=y_{p+1},\qq p=2N+1$.

 \no  ii) The spectrum of $H_b$ is given by
\[
\lb{spp}
 \s(H_b)=\bigcup_{k\in \Z_N}\s_k,\qq
 \s_k=\l_k([0,2\pi])=\ca [ \l_k^{-},\l_k^{+}],
  & k\ge0\\
 [\l_k^{+},\l_k^{-}], & k<0\ac,\qq k\in \Z_{N}=\{-N,...,N\},
\]
where $\l_{-N}(t)\le \l_{-N+1}(t)\le ...\le \l_{N}(t)$ are
eigenvalues of $J_t, t\in [0,2\pi)$. Moreover, $\l_{n}(\cdot)$ is
real analytic for any $t$, where $a_{2k-1}\ne0$, $k\in \N_N$.
\end{theorem}

{\no\bf Remark.} 1) Below we
 will sometimes write $\l_k(t,b,v)$, $\s_k(b,v)$, $J_t(b,v),..$
  instead of $\l_k(t)$, $\s_k$,
$J_t, ..$, when several potentials $v$ or magnetic fields are being
dealt with. The operator $J_t=J_t(b,v)$ is the symmetric $p\ts p$
matrix given by
\[
\lb{Jk1}
 J_t(b,v)=\left(\begin{array}{ccccccc}
                                    v_1 & a_1 & 0 & 0 &. .. &..& 0 \\
                                    a_1 & v_2 & 1 & 0 &.. & ..&0\\
                                    0 & 1 & v_3 & a_3 &.. & ..& 0\\
                                    0 & 0& a_3 & v_4 & 1 &.. &0 \\
                                    .. & .. & .. & .. & ..&..&..\\
                                    .. & .. & .. & .. & ..&..&1\\
                                    0 & .. &..& 0 & 0& 1 & v_p
                 \end{array}\right)
                 = J_t(b,0)+\diag(v_n)_1^p, \qqq v=(v_n)_1^p\in \R^p.
\]
\no 2) We take the specific boundary of the ZGNR, see Fig. 2. If we
change the boundary of the ZGNR, then the corresponding operator
$J_t$ will be more complicated. In the last case even the exact
calculation of the spectrum for the unperturbed operator $\D_0$ can
be the problem.

\no 3) If $\l_k(t)=\const$ for all $t\in [0,2\pi]$, then $\s_k$ is a
flat band. Otherwise, $\s_k$ is a non-flat spectral band. In the
periodic spectral theory  the basic problem is to describe all flat
bands. Examples of the dispersive curves for $\l_k(\cdot)$ are given
in Fig. 5.

We recall the  result about the spectrum of $\s(\D_0)$ of the
Laplacian $\D_0$ with $v=0$ and $b=0$ from \cite{KKu3} (for more
details see \er{T2-1}-\er{T2-3})
\begin{multline}
\lb{T2-2a}
 \s(\D_0)=\s_{ac}(\D_0)\cup\s_{pp}(\D_0),\ \ \s_{pp}(\D_0)=\{0\},\ \
 \s_{ac}(\D_0)=[-\l^0,\l^0]\sm [-s_1,s_1],\\
 where\qq \l^0=(5+4c_1)^{1\/2},\qq c_{\a}=\cos\frac{\a\pi}{N+1},
 \qq s_{\a}=\sin\frac{\a\pi}{N+1}, \ \ \a\in\R.
\end{multline}
Recall that $\Z_{N}=\{-N,...,N\}$.
 Now we describe the spectrum of the magnetic operator $\D_b$.

\begin{theorem}
\label{T3} i) Let $b\in \T$. Then the spectrum of $\D_b$ is given by
\[\lb{T3-1}
 \s(\D_b)=\s_{ac}(\D_b)\cup \s_{pp}(\D_b),\qqq
\s_{ac}(\D_b)=\bigcup_{k\in\Z_N\sm\{0\}}\l_k([0,2\pi]),\qqq
\s_{pp}(\D_b)=\{0\},
\]
where $\s_0(\D_b)=\{0\}$ is a flat band and $\l_k(\cdot)$ is defined
in Theorem \ref{T1}.

ii) Let $b\in \T$ and $b\to 0$. Then
$\s_{ac}(\D_b)=[\m_1^-,\m_2^-]\cup [\m_2^+,\m_1^+]$ and $\m_n^{\pm}$
has asymptotics
\[\lb{T3-2}
 \m_1^{\pm}=\pm(5+4c_1)^{1\/2}+O(b^2),\ \
 \m_2^{\pm}=\pm s_1\mp\frac{3c_1\sqrt{4-c_1^2}}{2(N+1)s_1}b+O(b^2).
\]
\end{theorem}

 We describe the spectrum of $H_b$ for $b$ and $v=(v_n)_1^p\in\R^p$
 for which $\s_{pp}(H_b)\ne\es$.

\begin{theorem} \label{T2}
 Let $(b,v)\in\T\ts\R^p$. Then $\s_{pp}(H_b)\ne \es$ iff
$v_{2n+1}=v_1$ for any $n\in\N_{N}$. Moreover, if $\s_{pp}(H_b)\ne
\es$, then $\s_{pp}(H_b)=\s_0(H_b)=\{v_1\}$ and each $|\s_n(H_b)|>0$
and $\s_n(H_b)\ss (v_1,+\iy), \s_{-n}(H_b)\ss (-\iy,v_1)$ for all
$n\in \N_N$.
\end{theorem}
{\no\bf Remark.} (1) Theorem \ref{T2} generalizes results from
\cite{KKu3}, devoted to the operator $\D_0+V$.

(2) If $\s_{pp}(H_b)=\{v_1\}$, then a flat band $\s_0(H_b)=\{v_1\}$
lies in the gap of the continuous spectrum, see Fig 3 and 4.

\begin{figure}[h]
\lb{f004}
         \centering\includegraphics[clip]{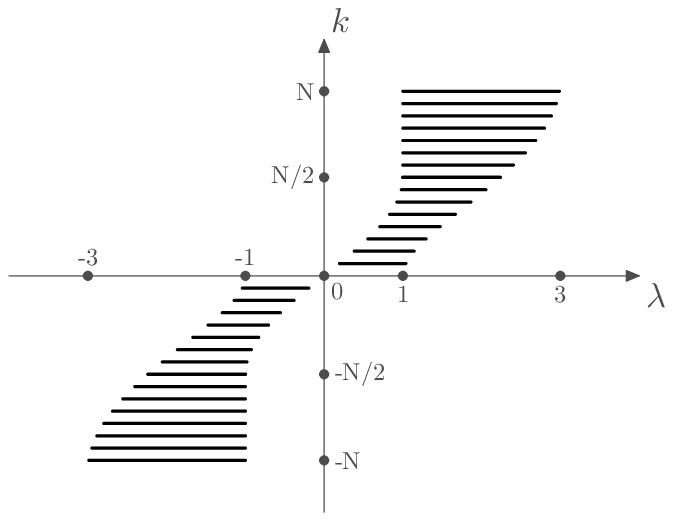}
         \begin{center}
         {\small Fig 3. Spectral bands $\s_k$ for the case $v,b=0$ and $N=15$. }
         \end{center}
\end{figure}

\begin{figure}[h]
\lb{f005}
         \centering\includegraphics[clip]{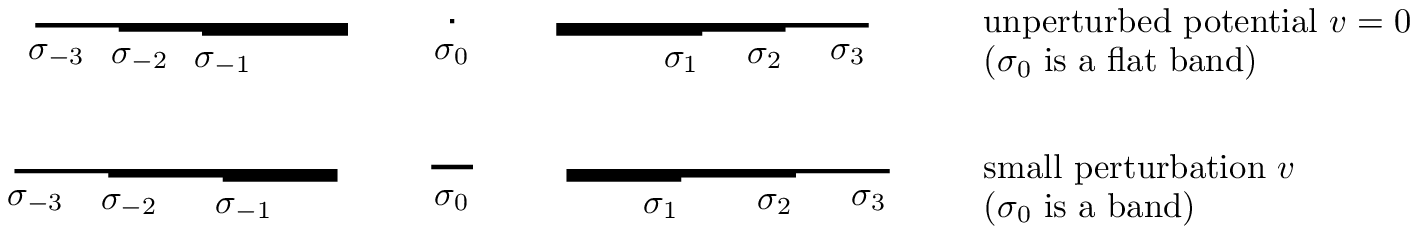}
         \begin{center}
         {\small Fig 4. Spectral bands for the case $N=3$. }
         \end{center}
\end{figure}
We will use notations $\l_k(t,b,v)\ev\l_k(t)$ for our spectral
curves $\l_k(t)$. The perturbation theory gives the standard
asymptotics
\[\lb{perc}
 \l_k(t,b,v)=\l_k(t,b,0)+O(\|v\|)\qqq as\ v\to0.
\]
In Theorem \ref{T4} we determine asymptotics of the eigenvalue
$\l_0(t,b,v)$ as $v\to0$ at fixed $b$, recall that $\l_0(t,b,0)=0$.
Below we will use this result to solve inverse spectral problem for
small potentials.

\begin{theorem}
\label{T4} Let $b\in\T$ and let $v\to0$. Then
\[\lb{T4-1}
 \l_0(t)=\sum_{n=0}^N v_{2k+1} \e_{k}(t)+O(\|v\|^2),
\]
\[\lb{T4-2}
 \eta_{k}(t)=\frac{\b^2_k(t)}{\sum_{s=0}^N\b_s^2(t)},\ \ \b_0=1,\ \
 \b_k=\prod_{j=1}^{2k}a_{j}(t),\ \ k\ge1,
\]
uniformly in $t\in[0,2\pi]$.
\end{theorem}

\begin{figure}[h]
\lb{f006}
         \centering\includegraphics[clip]{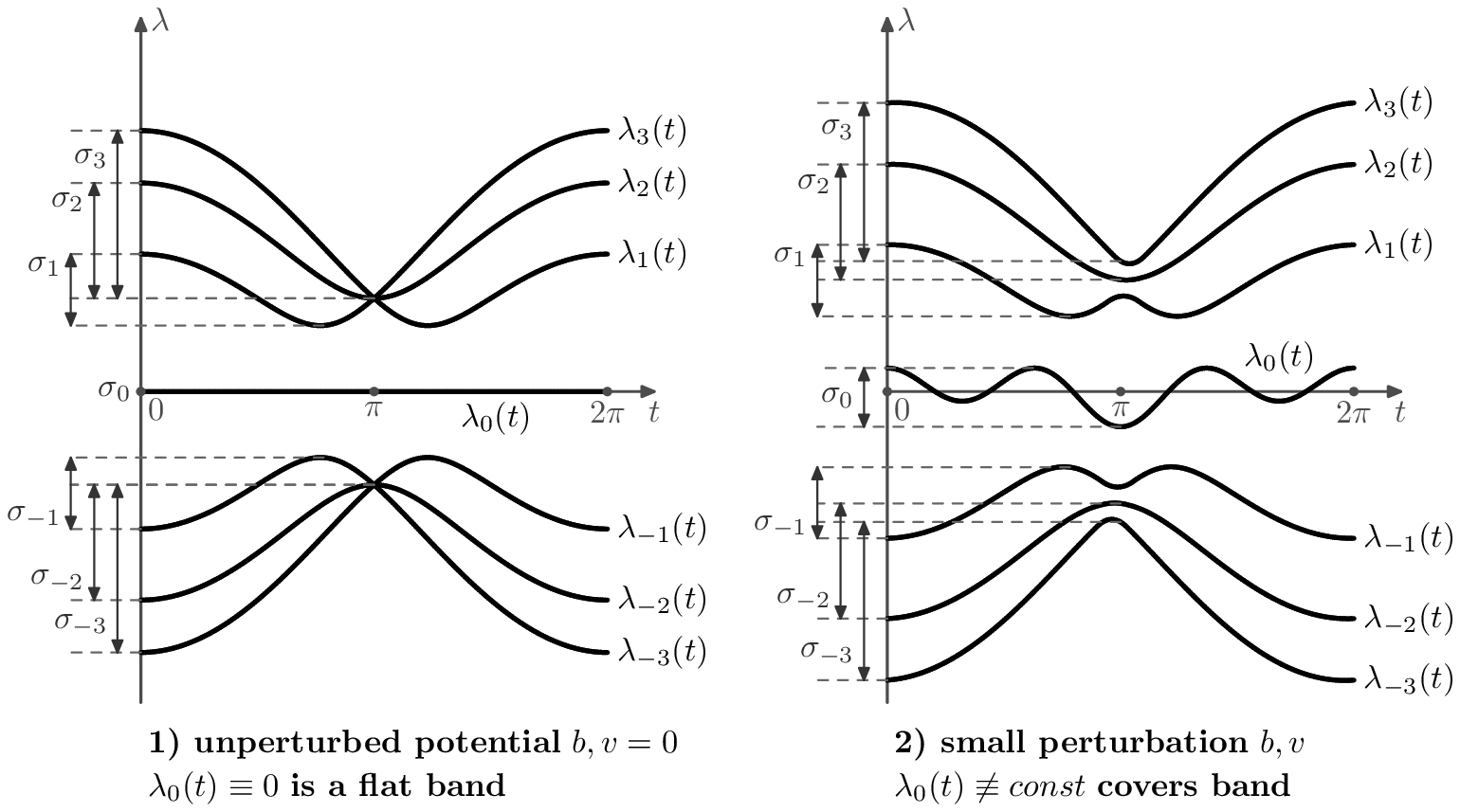}
         \begin{center}
         {\small Fig 5. Spectral curves $\l_k(t)$ and bands $\s_k$ for the case $N=3$. }
         \end{center}
\end{figure}

{\bf Strong electric fields.} We consider now the nanoribbon in
strong electric fields $\t V$ as the coupling constant $\t\to \iy$
and $V$ is fixed. Our operator has the form $H_b(\t)=\D_b+\t V$. We
determine asymptotics of the spectral bands $\s_{k}(\t
v)=[\l_{k}^+(\t),\l_{k}^{-}(\t)]$, $k\in\Z_N$ of the operator
$H_b(\t)$ as the coupling constant $\t\to \iy$.

\begin{figure}[h]
\lb{f007}
         \centering\includegraphics[clip]{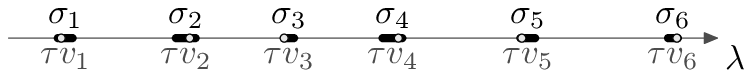}
         \begin{center}
         {\small Fig. 6. Spectral bands for the case of strong electric field }
         \end{center}
\end{figure}

\begin{theorem}
\lb{T5} Let $H_b(\t)=\D_b+\t V$, where $\t\to\iy$ and let the vector
$v=(v_k)_1^p\in \R^p, p=2N+1$ satisfy $v_1<..<v_p$. Then the
spectral bands $\s_{j}(\t v)=[\l_{j}^+(\t),\l_{j}^{-}(\t)],j\in
\Z_N$ have asymptotics:
\begin{multline}
\lb{T5.1}
 \l_{j}^{\pm s}(\t)=\t v_k-{\x_k^{\pm}+O(\t^{-1})\/\t},\qqq
 \x_k^{\pm}={r_{k}^{\pm}\/v_{k-1}- v_k}+{r_{k+1}^{\pm}\/v_{k+1} -
 v_k},\ \  s=(-1)^k,\qq k=j+N+1,\\
v_0=v_{p+1}=0,\qq  r_1^{\pm}=r_{p+1}^{\pm}=0,\ \ r_{2n+1}^{\pm}=1,\
\
 r_{2n}^-=0,\ \ r_{2n}^+=4,\ \ n\in\N_{N},
\end{multline}
\[
\lb{T5.3}
 |\s_{j}(\t v)|={4+O(\t^{-1})\/\t|v_{k-(-1)^k}-v_k|},\qq \ j\ne N,\qqq and \qqq
 |\s_{N}(\t v)|=O(\t^{-2}).
\]
Moreover, the operator $H_b(\t)$ has not a flat bands for $\t$ large
enough.
\end{theorem}

\no {\bf Remark.} 1) If $\t\to \iy$, then roughly speaking the
spectrum of operator $H_b(\t)$ consists of $p$ bands $\s_k(\t v),
k\in \N_{p}$ with lengths $|\s_k(\t v)|>0$, separated by $2N$ large
gaps. In this case all $2N$ gaps  are realized. In same sense it is
clear, since  we have the operator ${H_b\/\t}=V+\ve \D_b$ with small
coupling constant $\ve={1\/\t}$. The operator $V$ acting on
$\ell^2(\G)$ has only flat bands  (eigenvalues) $\{v_k\}, k\in
\N_p$. Under the small perturbation $\ve \D_b$ these flat bands
$\{v_k\}, k\ne N$ become the small spectral bands $\s_k(\t v)$ of
the continuous spectrum with lengths $>0$. It is important that only
Theorem \ref{T2} shows that $\s_N(\t v)$ is not a flat band, since
asymptotics \er{T5.3} are not sharp for $j=N$.

{\bf Inverse spectral problem for odd potentials.} Below we solve
inverse problems for sufficiently small "odd" potentials. Define the
space of odd potentials
$$
 \cV_{odd}=\{v=(v_k)_1^p\in\R^p:\ v_{2k}=0,\ k\in\N_N\}
$$
equipped with the standard norm $\|\cdot\|$ in $\R^p$. Let
$(b,v)\in\T\ts \cV_{odd}$ and let $0\le t_0<t_1<...<t_N\le \pi$. We
define the mapping $\L:\cV_{odd}\to\R^{N+1}$ by
$$
 v\mapsto\L(v)=(\l_0(t_j,b,v))_0^N,\qqq where \ \l_0(t_j,b,v)\in
 \s_0=\l_0([0,2\pi],b,v),
$$
and $\l_0$ is defined in Theorem \ref{T1}. Let
$B_{r}=\{v\in\cV_{odd}:\ \|v\|<r\}, r>0$ be the ball in $\cV_{odd}$.

\begin{theorem} \lb{T6}

i) Let $b\in\T$ be sufficiently small. Then the mapping $\L:B_{r}\to
\L(B_{r})$ is a real analytic bijection for some $r>0$, where
$\L(B_{r})\ss\R^{N+1}$ is an open domain and $0\in \L(B_{r})$.

ii) In the case of i) for any $\ve\in (0,r)$ there exists a
potential $v\in B_{\ve}$ and a potential $w\notin B_\ve$ such that
$\L(w)=\L(v)$, i.e. there is no a global injection.

iii) Let $b\in\T$ and let $t_k\in[0,2\pi)$, $k\in\N_{p}$ be a
sequence of a distinct numbers. Then the vector
$(\l_0(t_k,b,v))_1^{p}$ uniquely determines the potential $v\in
B_{r}$ for sufficiently small $r>0$.

\end{theorem}

%{\no\bf Remark.} We have shown that the spectral function
%$\l_0(t,b,V)\ev\l_0(t,V)$ uniquely determine "odd" potential.

%Assuming $f\in\ell^{\iy}(\G)$ in \er{H}, \er{V}, we may extend
%operator $H_b:\ell^2(\G)\to\ell^2(\G)$ to the operator
%$H_b:\ell^{\iy}(\G)\to\ell^{\iy}(\G)$.

%In the next Remark we use results from Theorem \ref{T6} and
%Proposition \ref{vec}.

\begin{figure}[h]
\lb{f006b}
         \centering\includegraphics[clip]{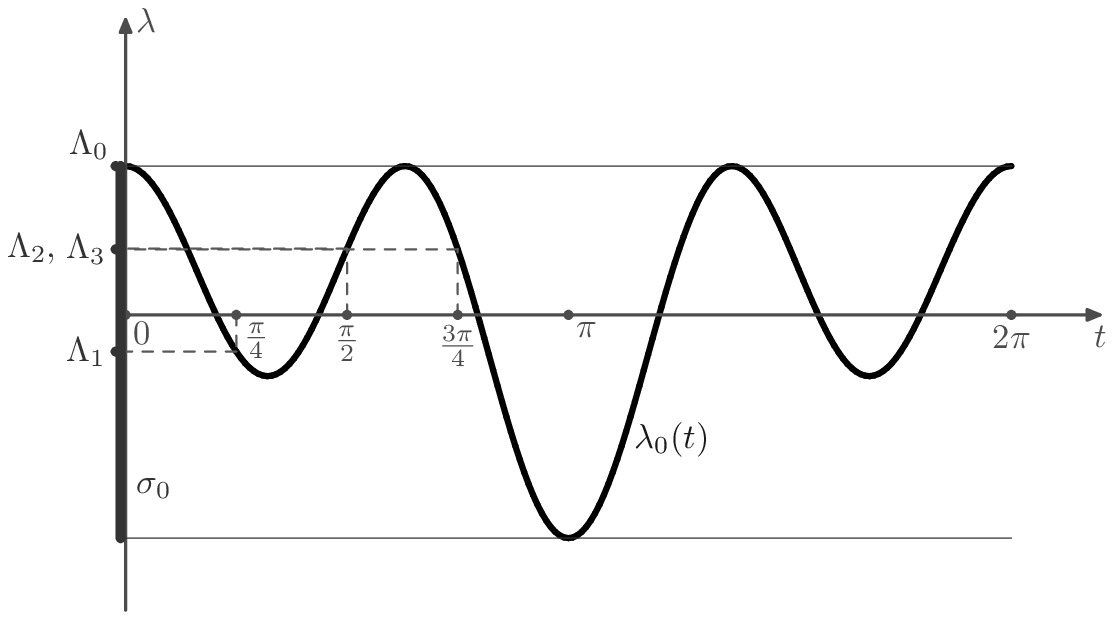}
         \begin{center}
         {\small Fig 6. The case $N=3$, small perturbation $b,v$. Spectral curve $\l_0(t)$ covers the band
         $\s_0$; $2(N+1)$-periodic eigenvalues $\L_k=\l_0(\frac{k\pi}{N+1})$ belong to $\s_0$.}
         \end{center}
\end{figure}

{\bf Remark.} Let $t_k=\frac{\pi k}{N+1}$, $k\in\N_N^0$. Then
$\l_0(t_k)$ are so-called $2(N+1)$-periodic eigenvalues for the
 Jacobi matrix $J_t$. Thus, if $b,
r$ are small enough, then by Theorem \ref{T6}.i), $2(N+1)$-periodic
eigenvalues $\l_0(t_k)\in \s_0$, $k\in\N_N^0$
   uniquely determine the odd potential $v$ (and
operator $V$). Note that we take only some part of $2(N+1)$-periodic
eigenvalues to determine the potential $V$, since $\l_j(t_k)$,
$j\ne0$ are $2(N+1)$-periodic eigenvalues too.

{\bf Inverse problems for monotonic electric  potentials without the
magnetic field.} We consider the case $b=0$. Define the mapping in
terms of antiperiodic eigenvalues
$$
\P: v\to  \P(v)=(\l_k(\pi,0,v))_{-N}^{N}.
$$
Here $\l_k(\pi,0,v)$ is a antiperiodic eigenvalues for Jacobi
operator $J_t$, since $t=\pi$. Define the bounded set of monotonic
potentials  by
$$
 \hat\cV_{\a}=\{v\in\R^p:\ 0\le v_1<v_2<v_3<..<v_{2N+1}\le\a\},
$$
$$
 \check\cV_{\a}=\{v\in\R^p:\ \a\ge v_1>v_2>v_3>..>v_{2N+1}\ge0\}.
$$

\begin{theorem}\lb{T7}\it
i) The mapping $\P:\hat\cV_{\a}\to\R^{N+1}$ is an injection iff
$\a\in[0,1]$.

ii) The mapping $\P:\check\cV_\a\to\R^{N+1}$ is an injection iff
$\a\in[0,1]$.

iii) In these both cases there exist the algorithm of recovering of
$v$ in terms of $\P(v)=(\l_k(\pi,0,v))_{-N}^{N}$.

iv) Spectral curves $\l_k$, $k\in\Z_N$ satisfy
$\l_k(\pi+t,0,v)=\l_k(\pi-t,0,v)$ for any $(t,v)\in[0,\pi]\ts\R^p$.
If $v\in\hat\cV_1$ or $v\in\check\cV_1$ then all $\l_k(\pi,0,v)$,
$k\in\Z_N$ are distinct numbers.
\end{theorem}

{\bf Remark.} 1) This Theorem shows that the set of all antiperiodic
eigenvalues uniquely determine the increasing (or decreasing)
bounded potential. Moreover, in the proof of Theorem \ref{T7} we
show how to determine the monotonic potential $v$ by antiperiodic
eigenvalues $\P(v)=(\l_k(\pi,0,v))_{-N}^{N}$.

2) Theorem \ref{T7}.iv) shows us that $\l_k$ reach the maximum or
minimum at the point $t=\pi$. Then for sufficiently small monotonic
potentials $v$ the spectrum change multiplicity at any of
anti-periodic eigenvalues $\l_k(\pi,0,v)$.

3) Note that in the paper [KKu2] we proved that some finite set of
$k$-periodic eigenvalues uniquely determines the spectrum (including
multiplicity) of a periodic matrix-valued Jacobi operator. In the
present paper Theorem \ref{T6}.i gives the stronger result  and
shows that some of $N+1$-periodic eigenvalues uniquely determine the
operator $H_b$ for potentials $v$  small enough.

In the proof of our theorems we determine various asymptotics for
periodic Jacobi operators with specific coefficients given by
\er{Jk}. Note that there exist a lot of papers devoted to
asymptotics and estimates for periodic Jacobi operators, see e.g.
\cite{KKu1}, \cite{KKr}, \cite{La}, \cite{vMou}.

We present the plan of our paper. In Section 4 we prove Theorem
\ref{T1} with a technical proof, where  we use arguments from
\cite{KKu3}. In Sect. 2 we describe the properties $\D_b$ and prove
Theorem \ref{T3}.  In Sect. 3 we prove Theorems \ref{T2} -\ref{T7}.

\section{Properties of the magnetic operator $\D_b$}
\setcounter{equation}{0}

Theorem 1.1 is proved in Section 4. In this section we describe the
spectral properties of $\D_b$.  Recall that Theorem 1.1 gives that
the operator $H_b=\D_b+V$ is unitarely equivalent to the operator
$\int_{[0,2\pi)}^{\os}J_t{dt\/2\pi }$, where $J_t\ev J_t(b,v)$  is a
$p\ts p$  Jacobi operator. Introduce a Jacobi operator $\wt J_t\ev
\wt J_t(b,v)$ acting in $\C^p$ and given by
\[
 \lb{i1100}
 (\wt J_t y)_n=\wt a_{n-1}y_{n-1}+\wt a_{n}y_{n+1}+v_ny_n,\qq
 \wt a_{n}=\ca  1 & {\rm even}\  n \\
2\cos(\frac t2-\frac{3n-2}2 b)& {\rm odd}  \ n\ac,
\]
and recall that $J_t\ev J_t(b,v)$ is given by
\[
 \lb{Jk2}
 (J_t y)_n=a_{n-1}y_{n-1}+a_{n}y_{n+1}+v_ny_n,\qq
 a_{n}=|\wt a_n|,
\]
where $y=(y_n)_1^{p}\in\C^{p}$ and $y_0=0=y_{p+1}, p=2N+1$. Note
that the matrices $\wt J_t=\wt J_t(b,v)$ and $J_t=J_t(b,v)$ are
unitarely equivalent. Then the eigenvalues $\l_k(t)=\l_k(t,b,v),
k\in\Z_N$ of $J_t$ (see Theorem \ref{T1}.ii) are eigenvalues of $\wt
J_t$ . Sometimes below we will use the matrix $\wt J_t$ instead of
$J_t$ to study the functions $\l_k(t)$ and the bands $\s_k$, defined
in Theorem \ref{T1}.ii).

{\bf The Laplacian $\D_0$.} We describe the spectrum $\s(\D_0)$ of
the unperturbed operator $\D_0$ at $V=0$ and $b=0$. Let
$\l_{-N}^0(t)\le \l_{-N+1}^0(t)\le ...\le \l_{N}^0(t)$ be
eigenvalues of the Jacobi matrix $J_t^0=J_t(0,0)$ corresponding to
$\D_0$ (see Theorem \ref{T1}.i). We recall the results about the
spectrum of $\s(\D_0)$ from \cite{KKu3}:
\[
\lb{T2-1}
 \l_k^0(t)=(a^2-2c_ka+1)^{1\/2}{k\/|k|},\ \
 1\le |k|\le N,\ \ \l_{0}^0(t)=0,\ \ t\in[0,2\pi],\ \
 a=2|\cos\frac t2|,
\]
\[
\lb{T2-2}
 \s^0=\s_{ac}^0\cup \s_{pp}^0=
 \bigcup_{-N}^{N} \s_k^0,\ \ \s_{pp}^0=\s_0^0=\{0\},\ \
 \s_{ac}^0=[-\l_N^0(0),\l_N^0(0)]\sm [-s_1,s_1],
\]
\[
\lb{T2-3} \s_{-k}^0=-\s_k^0,   \qq
\s_k^0=[\l_k^{0,-},\l_k^{0,+}]=\ca [\l_k^0(t_k^0), \l_k^0(0)],  &  if \ \ c_k<0\\
            [\l_k^0(\pi),\l_k^0(0)], &  if \ c_k\ge0 \ac,\ \ all \
            k=1,..,N,
\]
\[
\lb{T2-3a} \s_k^0=[\l_k^{0,-},\l_k^{0,+}]=-\s_{-k}^0=\ca
[\l_k^0(t_k^0),\l_k^0(0)]=[s_k, (5-4c_k)^{1\/2}],  &  if \ \ c_k\ge0 \\
            [\l_k^0(\pi),\l_k^0(0)]=[1,(5-4c_k)^{1\/2}], &  if \  c_k<0\ac,
\]
for all $k=1,..,N$, where $\l_k^0(t_k^0)=s_k$,
$2\cos\frac{t_k^0}2=c_k$, and $\l_k^0(\pi)=1$,
$\l_N^0(0)=(5+4c_1)^{1\/2}$.

 We use notation $\wt
J_t(b,v)$, $J_t(b,v)$ and $\l_k(t,b,v)$ for matrices $\wt J_t$,
$J_t$ (defined in \er{i1100}, \er{Jk}) and its eigenvalues.

For $k\ne0$ introduce
\begin{multline}
 \lb{i1}
 \F_k^0(t)=(\vp_{k,n}^0(t))_{n=1}^p,\ \
 \vp_{k,2n}^0(t)=\frac{(-1)^{n+1}s_{nk}}{\sqrt{N+1}},\ \
 \vp_{k,2n}^0(t+\pi)=\vp_{N+1-k,2n}^0(t),\ t\in (0,\pi),\\
 \vp_{k,2n-1}^0(t)=\frac{(-1)^{n+1}(2s_{nk}\cos\frac t2-s_{(n-1)k})}
 {\sqrt{N+1}\l_k^0(t)},\ \ \vp_{k,2n-1}^0(t+\pi)=\vp_{N+1-k,2n-1}^0(t),\
 t\in(0,\pi).
\end{multline}
Below we use notation $\langle \cdot,\cdot \rangle$ for the standard
scalar product in $\C^n$. We consider the case $V=0$ and $b\to0$.

\begin{lemma} \lb{LL1}
i) The following identities hold true
\[
 \lb{LLi1}
 J_t^0\F_k^0(t)=\l_k^0(t)\F_k^0(t),\ \ \|\F_k^0(t)\|=1,\ \
 k\ne0,\qqq
 t\in [0,2\pi],
\]
i.e. $\F_k^0(t)$, $k\ne0$ are orthonormal eigenvectors for $J_t^0$.

ii) \lb{LL3} For any $t\in[0,2\pi]\sm\{\pi\}$ the following identity
holds true
\[ \lb{i2}
 \pa_b(\l_k)(t,0,0)=\frac{\sin\frac t2}{(N+1)\l_k^0(t)}
 \lt((N+1)(3N+1)\lt(2\cos\frac t2+c_k\sign(t-\pi)\rt)-6\cos\frac
 t2\rt).
\]
\end{lemma}
{\no\bf Proof.}  Direct calculations gives i).

ii) Using the perturbation theory, we obtain
$$
 (\l_k)'_b(t,0,0)=\langle(\wt J_t)'_b(0,0)\F_k^0(t),\F_k^0(t)\rangle=
 2\sin\frac
 t2\sum_{n=1}^N(6n-5)\vp_{k,2n-1}^0(t)\vp_{k,2n}^0(t)
$$
$$
 =\frac{2\sin\frac t2}{(N+1)\l_k^0(t)}
 \sum_{n=1}^N(6n-5)(2s_{nk}^2\cos\frac t2-s_{nk}s_{(n-1)k})
$$
$$
 =\frac{2\sin\frac t2}{(N+1)\l_k^0(t)}
 \sum_{n=1}^N(12ns_{nk}^2\cos\frac t2 -
 10s_{nk}^2\cos\frac t2 -6ns_{nk}s_{(n-1)k}+5s_{nk}s_{(n-1)k})=
$$
$$
 =\frac{2\sin\frac t2}{(N+1)\l_k^0(t)}
 \lt((N+1)(3N-2)\cos\frac t2-(N+1)(3N+1)\frac{c_k}2\rt),
$$
where we use identities
$$
 \sum_{n=0}^N c_{2nk+\a}=\sum_{n=0}^N s_{2nk+\a}=0,\qqq \ 1\le|k|\le
 N,\ \ \a\in\R,
$$
$$
 \sum_{n=1}^N s_{nk}^2=\sum_{n=1}^N\frac{1-c_{2nk}}2=\frac{N+1}2,
$$
$$
 \sum_{n=1}^N
 s_{nk}s_{(n-1)k}=\frac12\sum_{n=1}^N(c_k-c_{(2n-1)k})
 =\frac12(Nc_k+c_{-k})=\frac{N+1}2c_k,
$$
$$
 \sum_{n=1}^Nnc_{2kn}=-\frac{N+1}2,\qqq
 \sum_{n=1}^Nns_{2kn}=-\frac{N+1}2\frac{c_k}{s_k},\ \ 1\le|k|\le N,
$$
$$
 \sum_{n=1}^Nn
 s_{nk}^2=\sum_{n=1}^Nn\frac{1-c_{2nk}}2=\frac{(N+1)^2}4,
$$
$$
 \sum_{n=1}^N
 n s_{nk}s_{(n-1)k}=\frac12\sum_{n=1}^Nn(c_k-c_{(2n-1)k})=
 \frac12\sum_{n=1}^Nn(c_k-c_{2nk}c_k-s_{2nk}s_k)=
$$
$$
 =\frac{(N+1)(N+2)}4c_k. \BBox
$$

\begin{lemma} \lb{LL2} Let $b, t\in[0,2\pi]$ and $v_{2n+1}=0$ for all
$n\in\N_N^0$. Then $\l=0$ is a simple eigenvalue of $J_t$.
\end{lemma}
\no {\bf Proof.} Define  a Jacobi matrix $A$ by
$$
 A=\left(\begin{array}{ccccc}
                                    0 & a_1 & 0 & .. & 0 \\
                                    a_1 & v_2 & a_2 & .. & 0\\
                                    0 & a_2 & 0 & .. & 0\\
                                    .. & .. & .. & .. & ..\\
                                    0 & .. & 0 & a_{2n} & 0
                 \end{array}\right),\qqq a_{2k}\ne0, for\ all \
                 k\in\N_n,
$$
for some $n\in\N$.
 Then $\l=0$ is a simple  eigenvalue of $A$ and $x=(x_k)_1^{2n+1}$ is a
corresponding eigenvector, where $x_1=1$,
$x_{2k+1}=(-1)^k\frac{a_1..a_{2k-1}}{a_2..a_{2k}}$, $x_{2k}=0$,
$k\in\N_n$. \BBox

{\bf Proof of Theorem \ref{T3}.} i) Lemma \ref{LL2} gives that
$\{0\}$ is a flat band. Then Theorem \ref{T2}. i) and the Theorem
\ref{T1} give the proof of i).

%The fact, that the value $0$ is a flat band follows from Lemma
%\ref{LL2}.

ii)  Relations \er{T2-1}-\er{T2-3} imply
\[
 \lb{t0000}
 \m_2^+(0)=\min_{t\in[0,2\pi]}\l_1^0(t)=\l_1^0(z_1^0)=\l_1^0(z_2^0)=s_1,\ \ \
 s_1<\l_1^0(t),\ \ t\in[0,2\pi]\sm \{z_1^0,z_2^0\},
\]
where $z_1^0,z_2^0\in(0,2\pi)$, $z_1^0\ne z_2^0$ such that
$2|\cos\frac{z_j^0}2|=c_1$, $j=1,2$. We denote
$\l_1^0(t,b)\ev\l_1(t,b,0)$, Perturbation Theory gives us that
$\l_1^0(t,b)$ is an analytic function in some neighbourhood of
points $t=z_1^0,z_2^0$ and $b=0$, because $\l_1^0(z_1^0,0)$ and
$\l_1(z_2^0,0)$ are simple eigenvalues of analytic matrix function
$\wt J_t(b,0)$. Using
\[\lb{loc2}
 \pa_t\l_1^0(z_j^0,0)=0,\ \
 \pa_t^2\l_1^0(z_j^0,0)>0,\ \ j=1,2,
\]
and the Implicit Function Theorem we deduce that there exist only
two functions $z_1(b),z_2(b)$ (analytic in $|b|\le \d$ for some
$\d>0$) such that
\[\lb{loc1}
 z_j(0)=z_j^0,\ \ \pa_t\l_1^0(z_j(b),b)=0,\ \
 \pa_t^2\l_1^0(z_j(b),b)>0,\ \ j=1,2
\]
for all $|b|\le \d$. Identities \er{t0000}-\er{loc1} yield that
$z_j(b)$ are points of local minimum of function $\l_1^0(t,b)$ and
one of these is a point of global minimum, that is
\[\lb{as03}
 \m_2^+(b)=\min_{t\in[0,2\pi]}\l_1^0(t,b)=\min_{j\in\N_2}\{\l_1^0(z_j(b),b)\},
\]
\[\lb{as02}
 \l_1^0(z_j(b),b)=\l_1^0(z_j^0)+b(\pa_b\l_1^0(z_j(b),b)|_{b=0})+O(b^2)=
 s_1+b(\pa_b\l_1^0(z_j(b),b)|_{b=0})+O(b^2),
\]
\[\lb{as01}
 \pa_b\l_1^0(z_j(b),b)|_{b=0}=
 \pa_t\l_1^0(z_j^0,0)\pa_b z_j(b)+\pa_b\l_1^0(z_j^0,0)
 =\pa_b\l_1^0(z_j^0,0),\ \ j=1,2.
\]
Then, substituting $\pa_b\l_1^0(z_j^0,0)=\pa_b\l_1(z_j^0,0,0)$ from
\er{i2} to \er{as01}, substituting \er{as01} to \er{as02} and using
\er{as03} we obtain asymptotics for $\m_2^+$. Other cases are proved
similar. \BBox

\section{Proof of the main results}
\setcounter{equation}{0}

%{\bf Proof of Theorem \ref{T2} ii)} Using famous results from
%perturbation theory we obtain $\l_k(t,b,v)=\l_k(t,0,0)+o(1)$ as
%$b\to0$, $v\to0$. And using \er{T2-1}-\er{T2-3}, we see that
%$\l_k(t,0,0)\ne\const$, $k\ne0$. Then $\l_k(t,b,v)\ne\const$ for
%sufficiently small $b$ and $v$. Then, using Theorem \ref{T1}, we
%obtain $\s_k$, $k\ne0$ are absolutely continuous spectral bands,
%since $\l_k(t,b,v)$ is analytic function on interval $t\in[0,2\pi]$
%without some finite set of points. \BBox

In order to prove Theorem \ref{T2} we need two Lemmas.

Recall that $v=(v_k)_1^p$, $p=2N+1$ and $a_{2k}=1$, $k\in\N_N$ and
$\wt a_{2k+1}=\wt a_{2k+1}(t,b)=2\cos(\frac t2-\frac{6k+1}2 b)$,
$k\in\N^0_N$. These functions are analytic for $t\in\C$, below we
will consider some cases when $t\in\C$ not only $t\in[0,2\pi]$ as we
did above. Introduce the matrices
\[
 \lb{bz1}
 T_k(t,b,v)=\ma -1 & -v_{2k-1}\\ v_{2k} & v_{2k}v_{2k-1}-\wt a_{2k-1}^2
 \am,\ \ k\in\N_{N+1},\ \ v_{2N+2}=0.
\]
Introduce the functions $u_k\ev u_k(t,b,v)$, $k\in\N^0_{N+1}$ by
\[ \lb{bz2}
 \lt(\begin{array}{c}u_{0}\\ u_{1}\end{array}\rt)=
 \lt(\begin{array}{c}0\\1\end{array}\rt),\ \ \
 \lt(\begin{array}{c}u_{2k}\\u_{2k+1}\end{array}\rt)=T_k
 \lt(\begin{array}{c}u_{2k-2}\\u_{2k-1}\end{array}\rt),\ \
 k\in\N_{N+1}.
\]
Note that $u_k$ is an analytic functions for $t\in\C$. Below we need
a simple fact.

\begin{lemma}\lb{Lz1}
Let $t,b\in [0,2\pi]$ and let $v\in \R^p$. Suppose $\wt
a_{2k-1}(t,b)\ne0$ for all $k=1,..,N $. Suppose $\l=0$ is eigenvalue
for $\wt J_t(b,v)$.
%(see Theorem \ref{T1} and eq.
%\er{i11}).
Then $u_{2N+2}(t,b,v)=0$.
\end{lemma}
\no {\bf Proof.} Let $\wt J_t(b,v)F=\l F=0$, where $F=(f_n)_1^p\ne0$
is the corresponding eigenvector. Then the definition  \er{Jk2}
implies
$$
 v_1f_1+a_2f_2=0,\ \ \
 \wt a_{2k-1}f_{2k-1}+v_{2k}f_{2k}+\wt a_{2k}f_{2k+1}=0,\ k\in\N_N,\ \
 \wt a_{2N}f_{2N}+v_{2N+1}f_{2N+1}=0,
$$
which yields
%(see \er{bz1}, \er{bz2} and we use $a_{2k}=1$,
%$k\in\N_N$)
$$
 f_1=f_1u_1,\ \ \ \lt(\begin{array}{c}f_{2k}\\f_{2k+1}\end{array}\rt)
 =\frac1{\wt a_{2k-1}}T_k
 \lt(\begin{array}{c}f_{2k-2}\\f_{2k-1}\end{array}\rt)=
{f_1\/\x_k}
 \lt(\begin{array}{c}u_{2k}\\u_{2k+1}\end{array}\rt),\ \
 k\in\N_{N+1},
$$
$\x_k=\prod_{n=1}^k \wt a_{2n-1}$   which yields $f_1\ne0$, since
$F\ne0$. Then
$$
 0=\wt a_{2N}f_{2N}+v_{2N+1}f_{2N+1}=
 {\wt a_{2N}u_{2N}+v_{2N+1}u_{2N+1}\/f_1 \x_N}=-
 {u_{2N+2}\/f_1 \x_N},
$$
which yields $u_{2N+2}=0$. \BBox

\begin{lemma}\lb{Lz2}
Let $b\in \T, v\in \R^p$. The following asymptotics hold true
\begin{multline}\lb{bz3}
 |u_{2k}(it)|=|v_{2k-1}|e^{2(k-1)t}+O\lt(e^{2(k-2)t}\rt),\ \ \
 t\to+\iy,\ \ \ i=\sqrt{-1},\ \ \ k\in\N_{N+1},\\
 |u_{2k+1}(it)|=e^{2kt}+O\lt(e^{2(k-1)t}\rt),\ \ \
 t\to+\iy,\ \ \ i=\sqrt{-1},\ \ \ k\in\N_{N}.
\end{multline}
\end{lemma}
\no {\bf Proof.} Note that
\[\lb{loc4}
 \wt
 a_{k}(it)=e^t+O(1),\ \ t\to+\iy,\ \ k\in\N_p.
\]
For $k=1$ identities \er{loc4} and \er{bz1}-\er{bz2} yield
\[\lb{loc5}
 \lt(\begin{array}{c}u_{2k}(it)\\ u_{2k+1}(it)\end{array}\rt)=
 (-1)^k\lt(\begin{array}{c} v_{2k-1}e^{(k-1)t}\\
 e^{kt}\end{array}\rt)+\lt(\begin{array}{c}O(e^{(k-2)t})\\
 O(e^{(k-1)t})\end{array}\rt),\ \ t\to+\iy.
\]
Substituting \er{loc5} for $k=1$ into \er{bz2} and using \er{bz1},
\er{loc4} we deduce that \er{loc5} is true for $k=2$. Repeating this
procedure we deduce that \er{loc5} is true for any $k\in\N_p$, which
yield asimptotic \er{bz3}. \BBox

%
%
%\begin{lemma}\lb{Lz2a}
%Let $b\in \T, v\in \R^p$. The following identities hold true:
%\begin{multline}\lb{bz3a}
% u_{2k}(t)\ev u_{2k}(t,b,v)=
%(-1)^{k} v_{2k-1}\prod_{j=0}^{k-1}\wt a_{2j-1}^2+r_{2k-1},\ \ \ k\in\N_{N+1},\\
% a_{-1}=1,\ \ r_{1}=0,\ \ r_{2k-1}(t)\ev r_{2k-1}(t,b,v)=O(e^{2(k-2)|t|}),
% \ \ it\to +\iy,
% i=\sqrt{-1},
%\end{multline}
%\begin{multline}\lb{bz4a}
% u_{2k+1}(t)\ev u_{2k+1}(t,b,v)=(-1)^k\prod_{j=1}^k
% \wt a_{2j-1}^2+\wt r_{2k+1}(t),\ \ \ k\in\N_N,\\ \wt r_{2k+1}(t)\ev\wt r_{2k+1}(t,b,v)
% =O(e^{2(k-1)|t|}),\ \ it\to  +\iy,\ \ i=\sqrt{-1}.
%\end{multline}
%\end{lemma}
%\no {\bf Proof.} For $k=1$ these identities are clear. Using
%\er{bz1}, \er{bz2} and $a_{2n-1}(t)\ev a_{2n-1}(t,b)\sim
%\frac{e^{|t|}}2$, $it\to+\iy$ we obtain \er{bz3}, \er{bz4} for
%$k\ge1$, step by step. ?????\BBox

\no {\bf Proof of Theorem \ref{T2}} If $v_1=0$ then {\bf
Sufficiency} follows from Lemma \ref{LL2}. If $v_1\ne0$ then the
proof is similar, since we may consider new operator
$H_b=\D_b+(V-v_1I)$, where $I$ is an identity operator.

{\bf Necessity.} We have $\s_{pp}\ne\es$. Without lost of generality
we may assume $0\in\s_{pp}$ (in other case we shift spectrum by
adding some constant $c$ to any of diagonal components $v_k$). Then
$\l=0$ is an eigenvalue for matrix $\wt J_t(b,v)$ for $t\in I_1$,
for some set $I_1\ss[0,2\pi]$ and $|I_1|>0$. There exists infinite
set $I_2\ss I_1$, satisfying $\prod_{n=1}^{N}\wt a_{2n-1}(t,b)\ne0$
for any $t\in I_2$. Then, by Lemma \ref{Lz1}, we have
$u_{2N+2}(t,b,v)=0$, $t\in I_2$. Then $u_{2N+2}(t,b,v)=0$ for any
$t\in\C$, since $u_{2N+2}$ is an analytic function by $t$. Using
\er{bz3} we deduce $|u_{2N+2}(it)|\sim |v_{2N+1}|e^{2Nt}$,
$t\to+\iy$, which yields $v_{2N+1}=0$. Then, using \er{bz1},
\er{bz2} and $v_{2N+1}=0$, we obtain
$u_{2N}(t,b,v)=-u_{2N+2}(t,b,v)=0$, $t\in\C$. Similarly to above, we
deduce $v_{2N-1}=0$ and so on.

By \er{T2-1} we know that
$\l_{-k}(t,0,0)<\l_0(t,0,0)\ev0<\l_k(t,0,0)$ for any
$(t,k)\in[0,2\pi]\ts\N_N$. Let some $b\in\T$ and let $v\in\R^p$ be
such that $v_1=v_3=..=v_{2N+1}=0$. We fix $t\in[0,2\pi]$ and let
variable $\t\in[0,1]$, perturbation theory gives us that $\l_k(t,\t
b,\t v)$, $k\in\Z_N$ are continuous functions depending on
$\t\in[0,1]$, since these are eigenvalues of analytic
matrix-function $\wt J_t(\t b,\t v)$. Note that $\l_0(\t b,\t v)=0$
for any $\t\in[0,1]$. Then inequalities $\l_{-k}(t,\t b,\t
v)<\l_0(t,\t b,\t v)\ev0<\l_k(t,\t b,\t v)$, $k\in\N_N$ remains for
any $\t\in[0,1]$, since by Lemma \ref{LL2} $\l_0=0$ is always simple
eigenvalue of $\wt J_t(\t b,\t v)$. Then by \er{spp}
$\s_{-k}(b,v)\ss(-\iy,0)$ and $\s_k(b,v)\ss(0,+\iy)$ for any
$k\in\N_N$. If $v_1=v_3=..=v_{2N+1}\ne0$ then we just shift spectrum
by adding $-v_1$ to any of diagonal components $v_k$ of the matrix
$\wt J_t$ and repeat the reasoning. \BBox

{\bf Proof of Theorem \ref{T4}}  If $V=0$, then Theorem \ref{T3}
gives that $\s_0=\{0\}$ is a flat band for $H^b$ . Also $\l_0=0$ is
a simple eigenvalue of matrix $J_t(b,0)$.
%Let $\l_0(t,b,v)$ be an
%eigenvalue of $\wt J_t(b,v)$, where $\l_0(t,b,0)=0$.
Using perturbation theory, for $v\to0$ we obtain
$$
 \l_0(t,b,v)=\l_0(t,b,0)+\sum_{n=1}^{2N+1}(\pa_{v_n}\l_0)(t,b,0)v_n+O(\|v\|^2),\
 \ uniformly\ for\ t\in[0,2\pi],
$$
$$
 (\pa_{v_n}\l_0)(t,b,0)=\langle(\pa_{v_n}J_t)(b,0)\F^b(t),\F^b(t)\rangle=(\vp_n^b)^2(t),
$$
where $\F^b(t)=(\vp_n^b(t))_1^p$ is a corresponding eigenvector for
an eigenvalue $\l_0(t,b,0)=0$, i.e. $J_t(b,0)\F^b(t)=0$ and
$\|\F^b(t)\|=1$. Using results from the proof of Lemma \ref{LL2}, we
obtain $\vp_{2n}^b(t)=0$, $\vp_{2n+1}^b(t)=(-1)^n\b_n(t)$ (recall
that $\b_n$ are defined in \er{T4-2}), since $J_t(b,v)$ has the form
\er{i11}.
%We have $\l_k(t,b,0)\ne0$, $k\ne0$, since by Lemma  \ref{LL2}
%$\l_0(t,b,0)=0$ is a simple eigenvalue of $J_t(b,0)$. Then, using
%\er{spp}, we obtain $\s_k\cap\s_0=\es$ for sufficiently small $V$.
\BBox

We start to prove our results about the inverse problems. It is
convenient to prove Theorem \ref{T6} after the proof of Theorem
\ref{T4}, since we use the same notations.

\no {\bf Proof of Theorem \ref{T6}.} We use notation from Theorem
\ref{T4}.

i) Let $b=0$. For any fix $t\in\R$ the function $\l_0(t,0,\cdot)$ is
real analytic for small arguments, since $\l_0(t,0,0)$ is a single
eigenvalue of analytic matrix-function $\wt J_t(0,0)$. Then for
small $v\in\cV_{odd}$ there exist Frechet derivative
$\pa_{v}\L(v):\cV_{odd}\to\R^{N+1}$ is a real analytic function.
Using \er{T4-1}, \er{T4-2} we obtain
$$
 D_0\ev\pa_{v}\L(0)=(\eta_{m}(t_n))_{n=0,m=0}^{N,N}\in\R^{(N+1)\ts(N+1)},
$$
where $\e_m$ are defined in \er{T4-2}. We need to show that $D_0$ is
an isomorphism, or we need to show that $D_0$ is an injection, since
$\dim\cV_{odd}=\dim\R^{N+1}$. Suppose $D_0v^0=0$ for some
$v^0=(v_j)_1^{p}\in\cV$. Then, using \er{T4-1}, \er{T4-2} and
\er{Jk} we deduce
\[\lb{i6001a}
 f(t_n)=0,\ n\in\N_{N}^0,\ where\
 f(t)=\sum_{m=0}^N a^m v_{2m+1},\ \ \ a=\left(2\cos\frac t2\right)^{2}.
\]
Then $f\ev0$ and $v_{2m+1}=0$, $m\in\N_{N}^0$, since this is a
polynomial of degree $N$ has $N+1$ zeroes. We have $D_0$ is an
isomorphism and $\L(\cdot)$ is a bijection for small arguments. If
$b\ne0$ is a sufficiently small value then $D_0$ is an isomorphism
too, since $D_0$ analytically depends on $b$, then $\L$ is a
bijection for small arguments.

ii) Define the potential $v^0=(v^0_n)_{1}^{p}\in\cV_{odd}$, where
$v_1^0=-1$, $v_{2n+1}^0=n$, $v_{2n}^0=0$, $n\in\N_N$. Let $s>0$ be
sufficiently small. Then using \er{Jk} we obtain
\[\lb{n001}
 J_t(b,s^{-1}v^0)=s^{-1}A(s),\ \
 A(s)=\diag(v^0)+sJ_t(b,0).
\]
Let $\wt\l_{-N}(s)\le..\le\wt\l_{N}(s)$ be eigenvalues of the matrix
$A(s)$. We have $A(0)=\diag(v^0)$ and then $\wt\l_{-N}(0)=-1$,
$\wt\l_{-N+n}(0)=0$, if $n\in\N_N$ and $\wt\l_{n}(0)=n$, $n\in\N_N$.

Consider the asymptotic for $\wt\l_{-N+n}(s)$, $n\in\N_N$ for
$s\to0$. We have $\wt e_n=(\d_{2n,j})_{j=1}^{p}$ (where $\d_{i,j}$
is a Kronecker symbol) are eigenvectors of $A(0)$, corresponding to
eigenvalues $\wt\l_{-N+n}(0)$. Using famous results of the
perturbation theory, we obtain that the numbers $\wt\l_{-N+n}'(0)$,
$n\in\N_N$ are eigenvalues of matrix $\wt A=((A'(0)\wt e_j,\wt
e_k))_{j,k=1}^N$. Easy calculations gives $\wt A=0$, then
$\wt\l_{-N+n}'(0)=0$, $n\in\N_N$, and then
\[\lb{n002}
 \wt\l_{-N+n}(s)=\wt\l_{-N+n}(0)+s\wt\l_{-N+n}'(0)+O(s^2)=O(s^2)
 \ as \ s\to0,\ \ n\in\N_N.
\]
Using the fact, that $\l_n(t,b,s^{-1}v^0)$, $n\in\Z_N$ are
eigenvalues of $J_t(b,sv^0)$ and using \er{n001}, \er{n002} we
obtain
\[\lb{n003}
 \l_0(t,b,s^{-1}v^0)=s^{-1}\wt\l_0(s)=O(s),\ \ s\to0.
\]
Then $\L(s^{-1}v^0)=O(s)$, $s\to0$. Using i) we obtain, that for any
sufficiently small $s>0$ there exists unique sufficiently small
potential $w^s\in\cV_{odd}$, which satisfy $\L(s^{-1}v^0)=\L(w^s)$,
and $w^s\ne s^{-1}v^0$, since $s^{-1}v^0$ has large norm.

iii) Define the mapping $\wt\L(v)=(\l_0(t_k,b,v))_0^{p}$. We use
similar arguments as in i). For any fix $(t,b)\in\R^2$ the function
$\wt\l_0(.)\ev\l_0(t,b,\cdot)$ is real analytic for small arguments,
since for small $v$ the number $\l_0(t,b,v)$ is a single eigenvalue
of analytic matrix-function $J_t(b,v)$. Then there exists Frechet
derivative $\pa_{v}\wt\L(v):\cV_{odd}\to\R^{2N+1}$, which is real
analytic function for sufficiently small $v\in\cV_{odd}$. Using
\er{T4-1}, \er{T4-2} we obtain
$$
 D_0\ev\pa_{v}\wt\L(0)=(\eta_{2m+1}(t_n))_{n=0,m=0}^{2N,N}\in\R^{(2N+1)\ts(N+1)}.
$$
We need to show that the matrix (operator) $D_0$ is an injection.
Suppose $D_0v^0=0$ for some $v^0=(v_j)_1^{2N+1}\in\cV_{odd}$. Then,
using \er{T4-1}, \er{T4-2}, we deduce
\[\lb{i6001}
 f(t_n)=0,\ n\in\N_{2N}^0,\ where\ f(t)=\sum_{m=0}^N \b_m^2(t)v_{2m+1}.
\]
Using \er{T4-2} and \er{Jk} we obtain
\[\lb{i6002}
 \b_m^2(t)=\sum_{j=0}^m (c_{j,m}\cos jt+s_{j,m}\sin jt),\ \
 c_{m,m}\ne0,\ \ m\in\N_N^0
\]
for some constants $c_{j,m},s_{j,m}$. Then using \er{i6001},
\er{i6002} we deduce that $f(t)\ev0$, since $f$ is a trigonometrical
polynomial of degree $N$ and it has $2N+1$ zeroes in the interval
$[0,2\pi]$. Also we deduce that $v_{2m+1}=0$, $m\in\N^0_N$, since by
\er{i6002} we have that $\b_m^2(t)$, $m\in\N^0_N$ are linearly
independent functions. Then $v^0=0$ and $D_0=\pa_{v}\wt\L(0)$ is an
injection. Also $\wt\L(\cdot)$ is an injection for small arguments,
since $\pa_v\wt\L(\cdot)$ is a real analytic function for small
arguments. \BBox

We consider the case of the strong electric fields.

\no {\bf Proof of Theorem \ref{T5}.} In the case $\D_b+\t V$ the
corresponding Jacobi operator depends on $\t$ and is given by
\[
\lb{jt} (\wt J_t(b,\t v)y)_n=\wt a_{n-1}y_{n-1}+\wt a_{n}y_{n+1}+\t
v_ny_n,\ y=(y_n)_{n\in \Z}\in \C^p.
\]
Using \er{i1100} we rewrite $\wt J_t(b,\t v)$, $\t\to \iy$ in the
form
$$
 \wt J_t(b,tv)=\t(\wt V+\t^{-1}\wt J_t(b,0))=\t(\wt V+\ve\wt J_t(b,0)),\
                                  \ \wt V=\diag(v_j)_{1}^{p},\ as \
                                   \ \ve={1\/\t}\to 0.
$$
Then the perturbation theory (see  Sect. XII, 1, \cite{RS})  for
$\wt V+\ve\wt J_t(b,0)$ gives
\begin{multline}
\lb{T5.b}
 {\l_{k-N-1}(t,b,\t v)}=\t(v_k-\a_k\ve^2+O(\ve^3)),\\
 \a_k=\sum_{j\in \N_p\sm\{k\}}{u_{k,j}^2\/v_j-v_k},\ \
 u_{k,j}=(e_j^0,\wt J_t(b,0)e_k^0),
\end{multline}
for any $k\in\N_p$, here  $\wt V e_j^0=v_je_j^0$ and the vector
$e_j^0=(\d_{j,n})_{n=1}^{p}\in \C^{p}$. The matrix $\wt
J_t(b,0)=\{u_{k,j}\}$ is given by \er{i1100}, where
\[ \lb{T5.b1}
 u_{k,j}\ev u_{k,j}(t,b)=\ca0, & |k-j|\ne1,\\
 2\cos(\frac t2-\frac{6k-5}2 b), & j=k+1,\ j\in2\N,\ \\
 2\cos(\frac t2-\frac{6j-5}2 b), & j=k-1,\ k\in2\N
 \\ 1, & other\ cases \ac\ .
\]
Then \er{T5.b}, \er{T5.b1}  yield \er{T5.1}, \er{T5.3}. \BBox

{\bf Proof of Theorem \ref{T7}}. i) Sufficiency. Sufficiently to
show that $\P$ is an injection for $\a=1$, since
$\hat\cV_{\b}\ss\hat\cV_{\a}$ for $\b\le\a$. Using \er{Jk} we get
$a_{2n-1}(\pi,0)=0$, $a_{2n}=1$, $n\in\N_N$, so direct calculations
give
\[
\lb{loc10}
 \bigcup_{k=-N}^{N}\{\l_k(\pi,0,v)\}=\s(J_{\pi}(0,v))=
 \{v_1\}\cup\bigcup_{k=1}^N\{\l:\
 (\l-v_{2k})(\l-v_{2k+1})-1=0\}.
\]
Below we use the notation $\l_k^{(0)}\ev\l_k(\pi,0,v)$. Firstly, let
$p_j(\l)=(\l-x_j)(\l-y_j)-1$, $j=1,2$ be two quadratic polynomials
for some $0\le x_1<y_1<x_2<y_2\le1$. Let $\m_j<\n_j$ be roots of
$p_j$, then it is not difficult to show that
\[\lb{loc11}
 \m_1<\m_2<0<1<\n_1<\n_2.
\]
Using this fact, monotonicity of $v$ and \er{loc10} we deduce that
\[\lb{loc12}
 \l_{-N}^{(0)}<..<\l_{-1}^{(0)}<0\le\l^{(0)}_0\le 1<\l^{(0)}_1<..<\l_N^{(0)},
\]
where $\l^{(0)}_0=v_1$ and $\l_{-N-1+k}^{(0)}$, $\l_k^{(0)}$ are
roots of the polynomial $(\l-v_{2k})(\l-v_{2k+1})-1$ for $k\in\N_N$.
Then $\P(v)$ uniquely determine monotonic potential $v$ and then
$\P:\cV_{\a}\to\R^p$ is an injection.

Necessity. Suppose that $\a>1$. Define two potentials
$v,w\in\hat\cV_{\a}$ by
\begin{multline}
\lb{loc13}
 0<v_1=2\ve<\frac14<v_2<..<v_{p-2}<\frac34<\\
  v_{p-1}=1+5\ve-\sqrt{2\ve+\ve^2}<v_{p}=1+5\ve+\sqrt{2\ve+\ve^2}<\min\lt\{\a,\frac54\rt\},
\end{multline}
\begin{multline}\lb{loc14}
 0<w_1=4\ve<\frac14<w_2=v_2<..<w_{p-2}=v_{p-2}<\frac34<\\
  w_{p-1}=1+4\ve-2\sqrt{\ve+\ve^2}<w_{p}=1+4\ve+2\sqrt{\ve+\ve^2}<\min\lt\{\a,\frac54\rt\},
\end{multline}
for some $\ve>0$ small enough. Direct calculations give us that
$\P(v)=\P(w)=(\l^{(0)}_k)_{-N}^N$, where $\l^{(0)}_{-N-1+k}$,
$\l^{(0)}_k$ are roots of the polynomial
$(\l-v_{2k})(\l-v_{2k+1})-1$ for $k\in\N_{N-1}$ and
$\l^{(0)}_{-1}=2\ve$, $\l^{(0)}_{0}=4\ve$, $\l^{(0)}_{N}=2+6\ve$.
Then $\P:\cV_{\a}\to\R^{2N+1}$ is not injection, since
$\P(v)=\P(w)$. The Proof of ii) is similar to the Proof of i).

iii) Firstly, $\l^{(0)}_0=v_1$ determine the first component of
potential. By the remark after \er{loc12} the components $v_{2k}$,
$v_{2k+1}$ are defined uniquely as a roots of polynomial
$(\l-\l_{-N-1+k}^{(0)})(\l-\l_k^{(0)})+1$.

iv) By \er{Jk} we have that $J_{\pi+t}(0,v)=J_{\pi-t}(0,v)$, which
yields symmetry $\l_k(\pi-t,0,v)=\l_k(\pi+t,0,v)$ and then $t=\pi$
is a point of local extremum for functions $\l_k$. If
$v\in\hat\cV_1$ (or $\check\cV_1$) then by \er{loc12} all components
$\l_k^{(0)}$ are distinct numbers. \BBox

\section{ Appendix, proof of Theorem 1.1}
\setcounter{equation}{0}

For the magnetic field $\mB=B(0,0,1)\in \R^3$ the corresponding
magnetic vector potential is given by
$$
\mA(x)={1\/2}[\mB,x]={B\/2}(-x_2,x_1,0),\qqq x=(x_1,x_2,x_3)\in\R^3.
$$

Define the coefficients $a_{\n,\r}(t)=(\mA(\vk_\n+t{\bf
e}_{\n,\r}),{\bf e}_{\n,\r})$, where $\n, \r\in\Z\ts\N_p$ and $t\in
[0,1]$. In Lemma \ref{maf00} we will show that $a_{\n,\r}(t)$ does
not depend on $t\in [0,1]$. Recall that the magnetic operator $\D_b$
is given by
\begin{multline}
 \lb{Hlemma1}
 (\D_bf)_{\o}=
 e^{ia_{\o, (n,2k)}}f_{n,2k}+e^{ia_{\o,(n-1,2k+2)}}f_{n-1,2k+2}+
 e^{ia_{\o, (n,2k+2)}}f_{n,2k+2},\qq
 k\in \N^0_{N},
 \\
 f_{n,0}=f_{n,2N+2}=0, \qqq \qqq b={B\sqrt3\/2},
 \\
 (\D_bf)_{\s}=e^{ia_{\s, (n,2k-1)}}f_{n,2k-1}+e^{ia_{\s,(n+1,2k-1)}}f_{n+1,2k-1}
 +e^{ia_{\s, (n,2k+1)}}f_{n,2k+1},
 \ \ k\in\N_N,
\end{multline}
where  $\o=(n,2k+1)$, $\s=(n,2k)\in \Z\ts \N_p$ and
$f=(f_{n,k})_{(n,k)\in \Z\ts \N_p}\in\ell^2(\G)$.

\begin{lemma}

\lb{maf00} Let a function $a_{\n,\r}(t)=(\mA(\vk_\n+t{\bf
e}_{\n,\r}),{\bf e}_{\n,\r})$, where $\n, \r\in\Z\ts\N_p$ and $t\in
[0,1]$. Denote $\o=(n,2k+1), \s=(n,2k)$, then
\[
\lb{maf1} a_{\o, (n,2k+2)}=b(n-k),\qq a_{\o,(n-1,2k+2)}=b(n+2k),\qq
a_{\o, (n,2k)}=-b(2n+k),
\]
\[
\lb{maf2} a_{\s, (n,2k-1)}=-b(n-k+1),\qq
a_{\s,(n+1,2k-1)}=-b(n+2k-1),\qq a_{\s, (n,2k+1)}=b(2n+k),
\]
where $\qq b={B\sqrt3\/2}$ and all $(t,n,k)\in [0,1]\ts\Z\ts \N_p$.
%\[
%\lb{maf2} ??????? a_{\n, (n,2k+2)}={B\sqrt3}(n-k),\qqq
%a_{\o,(n-1,2k+2)}={B\sqrt3}(n+2k),\qq a_{\o,
%(n,2k)}={B\sqrt3}(2n+k),
%\]

\end{lemma}
\no{\bf Proof.} Identity $\mA({\bf r})={B\/2}[{\bf e}_0,{\bf r}],
{\bf e}_0=(0,0,1), {\bf r}\in \R^3$ yields for any $ t\in[0,1]$
\[
\lb{ao1} a_{\o,\s}(t)={B\/2}([{\bf e}_0,\vk_\o+t{\bf
e}_{\o,\s}],{\bf e}_{\o,\s})= {B\/2}([{\bf e}_0,\vk_\o],{\bf
e}_{\o,\s})=a_{\o,\s}(0)=a_{\o,\s},
\]
where $\s=(n,2k+2), (n-1,2k+2), (n,2k)$.

Recall $\vk_\o=(\sqrt3 (2n+k),3k,0)$. If $\s=(n,2k+2)$, then ${\bf
e}_{\o,\s}=({\sqrt3\/2},{1\/2},0) $ and \er{ao1} yields
$$
a_{\o, (n,2k+2)}={B\/2}([{\bf e}_0,\vk_\o],{\bf e}_{\o,\s})=
{B\/2}\det
\ma {\sqrt3\/2} & {1\/2} & 0\\
0 & 0 & 1\\
   \sqrt3 (2n+k) & 3k & 0\am={\sqrt3
B\/2}(n-k).
$$
If $\s=(n-1,2k+2)$, then ${\bf e}_{\o,\s}=(-{\sqrt3\/2},{1\/2},0) $
and \er{ao1} yields
$$
a_{\o,(n-1,2k+2)}={B\/2}([{\bf e}_0,\vk_\o],{\bf e}_{\o,\s})=
{B\/2}\det
\ma  -{\sqrt3\/2} & {1\/2} & 0\\
0 & 0 & 1\\     \sqrt3 (2n+k) & 3k & 0\am={B\sqrt3\/2}(n+2k).
$$
If $\s=(n,2k)$, then ${\bf e}_{\o,\s}=(0,-1,0) $ and \er{ao1} yields
$$
a_{\o,(n-1,2k+2)}={B\/2}([{\bf e}_0,\vk_\o],{\bf e}_{\o,\s})=
{B\/2}\det
\ma  0 & -1 & 0\\
0 & 0 & 1\\     \sqrt3 (2n+k) & 3k & 0\am=-{B\sqrt3\/2}(2n+k).
$$
The proof of other cases is similar. \BBox

Substituting identities from Lemma \ref{maf00} into the magnetic
operator $\D_b$ is given by \er{Hlemma1} we obtain \er{H}.

{\bf Proof of Theorem \ref{T1}.} i) Recall that
\begin{multline}
\lb{H2}
 (\D_bf)_{n,2k+1}=
 e^{-ib(2n+k)}f_{n,2k}+e^{ib(n+2k)}f_{n-1,2k+2}+
 e^{ib(n-k)}f_{n,2k+2}+v_{2k+1}f_{n,2k+1},\
 k\in \N^0_{N},
 \\
 f_{n,0}=f_{n,2N+2}=0,
 \\
 (\D_bf)_{n,2k}=e^{-ib(n-k+1)}f_{n,2k-1}+e^{-ib(n+2k-1)}f_{n+1,2k-1}
 +e^{ib(2n+k)}f_{n,2k+1}+v_{2k}f_{n,2k},
 \ k\in\N_N,
\end{multline}
Define the unitary operators $U, S $ acting in $\ell^2(\Z)$ by
\[\lb{vec2}
U(h_n)_{n\in \Z}=(\t^nh_n)_{n\in \Z},\qqq   \t=e^{ib}, \qqq \qqq
Sh=(h_{n+1})_{n\in\Z}, \qqq h=(h_n)_{n\in\Z}
\]
 For each
$(f_{n,k})_{(n,k)\in \Z\ts \N_p}\in \ell^2(\G)$ we introduce the
function $\p_k=(f_{n,k})_{n\in\Z}\in\ell^2(\Z)$, $k\in\N_{p}, p=2N+1
$ and $\p=(\p_k)_{k\in\N_p}\in(\ell^2(\Z))^p$. Using \er{H2} and
$(Vf)_{n,k}=v_kf_{n,k}$ for any $(n,k)\in\Z\ts\N_p$, we obtain that
the operator $H:\ell^2(\G)\to\ell^2(\G)$ is unitarily equivalent to
the operator $K:\ell^2(\Z)^p\to\ell^2(\Z)^p$, given by
$$
 (K\p)_{2k+1}=\ol\t^{k}(U^*)^2\p_{2k}+
 (\ol\t^{k}+S^*\t^{2k+1})U\p_{2k+2}+v_{2k+1}\p_{2k+1},\qq \p_0=\p_{p+1}=0,
 \qq k\in \N^0_{N},
$$
$$
 (K\p)_{2k}=U^*(\t^{k-1}+S\ol\t^{2k-1})\p_{2k-1}+U^2\t^{k}\p_{2k+1},
 \ \ k\in\N_N.
$$
 We
rewrite $K$ in the matrix form by
\[
 \lb{Km}
 K(\p_k)_{1}^{p}=\ma v_1 & (1+\t^{} S^*)U & 0 & 0&.. & 0 \\
U^*(I+\ol\t^{} S) & v_2 & U^2\t & 0&.. & 0\\
0 & (U^*)^2\ol\t & v_3 & (\ol\t+\t^{3} S^*)U & ..&0\\
.. & .. & .. & .. & .& U^2\t^N \\
0 & .&. & 0 & (U^*)^2\ol\t^N & v_{p}\am
\ma \p_1\\ \p_2\\ \p_3\\..\\
\p_{p}\am,
\]
where $ \p_k\in\ell^2(\Z)$. Note  that $K^*=K$, since $S^*=S^{-1}$.
We rewrite $K$ in the matrix form by
\[
 \lb{Km1}
 K=\ma v_1 & A_1 & 0 & 0&.. & 0 \\
A_1^* & v_2 & A_2 & 0&.. & 0\\
0 & A_2^* & v_3 & A_3 & ..&0\\
.. & .. & .. & .. & .&A_{p-1}\\
0 & .&. & 0 & A_{p-1}^* & v_{p}\am,
\]
where $A_{2k}=\t^kU^2$, $k\in\N_N$ and
$A_{2k+1}=(\ol\t^k+\t^{2k+1}S^*)U$, $k\in\N^0_N$. Define the unitary
operator $\mU=\diag (u_k)_1^p$, where $u_{2k}=U^{3k-2}$, $k\in\N_N$
and $u_{2k+1}=U^{3k}$, $k\in\N^0_N$. Using $US^*=\t SU^*$ we obtain
\[\lb{vec3}
K_1=\mU K\mU^*=\ma v_1 & r_1 & 0 & 0&.. & 0 \\
r_1^* & v_2 & r_2 & 0&.. & 0\\
0 & r_2^* & v_3 & r_3 & ..&0\\
.. & .. & .. & .. & .&r_{p-1}\\
0 & .&. & 0 & r_{p-1}^* & v_{p}\am,
\]
where $r_{2k}=\t^k$, $k\in\N_N$ and
$r_{2k+1}=\t^{-k}(\t^{6k+1}S^*+1)$, $k\in\N^0_N$.
%$$
%r_1=(1+\t^{} S^*),\ r_2=\t,\qq r_3=(\ol\t+\t^{3} S^*) ,\qq
%r_4=\t^2,\qq r_5=(\ol\t^2+\t^{5} S^*) ,...
%$$
%
%????  KONEC  ????????????????????????

Introduce the unitary operator  $\F:\ell^2(\Z)^p\to
\int_{[0,2\pi)}^{\os}\mH_0 {dt\/2\pi},\ \ \mH_0=\C^p$, by $\F
(\p_k)_{1}^p=(\f \p_k)_{1}^p$, where  $\f:\ell^2(\Z)\to L^2(0,2\pi)$
is an unitary operator given by
$$
\f  h=\sum_{n\in\Z}h_n{e^{int}\/\sqrt{2\pi}}, \qqq h=
(h_n)_{n\in\Z}\in\ell^2(\Z),\  t\in[0,2\pi].
$$
Then we deduce that
$$
\F_pK_1\F_p^{-1}=\int_{[0,2\pi)}^{\os}\wt{\wt J_t}{dt\/2\pi},
$$
where the operator $\wt{\wt J_t}: \mH_0\to \mH_0$ has the matrix
given by
\[\lb{vec1}
\wt{\wt J_t}y= \left(\begin{array}{ccccc}
                                    v_1 & 1+\t e^{-it} & 0 & .. & 0 \\
                                    1+\ol\t e^{it} & v_2 & \t & .. & 0\\
                                    0 & \ol\t & v_3 & .. & 0\\
                                    .. & .. & .. & .. & ..\\
                                    0 & .. & 0 & \ol\t^N & v_{p}
                 \end{array}\right)
                 \left(\begin{array}{c}y_1\\ y_2\\ y_3\\..\\ y_{p}\end{array}\right),
                 \qq y=(y_k)_1^{p}\in \C^{p}.
\]
The matrix $\wt{\wt J_t}$ is unitarely equivalent to the matrix
$J_t$, given by
\[\lb{i11}
 J_t=\left(\begin{array}{ccccc}
                                    v_1 & a_1(t) & 0 & .. & 0 \\
                                    a_1(t) & v_2 & a_2(t) & .. & 0\\
                                    0 & a_2(t) & v_3 & .. & 0\\
                                    .. & .. & .. & .. & ..\\
                                    0 & .. & 0 & a_{p-1}(t) & v_{p}
                 \end{array}\right),
\]
where $a_{2k}(t)=1$, $k\in\N_N$ and $a_{2k-1}(t)=2|\cos(\frac
t2-\frac{6k-5}2 b)|$, $k\in\N_N$.

Thus we deduce that the operator $H$ is unitarily equivalent to the
operator $\int_{[0,2\pi)}^{\os}J_t{dt\/2\pi}$.

ii) Let $\l_k(t), k\in \Z_N$ be eigenvalues of $J_t$ satisfy
$\l_{-N}(t)\le \l_{-N+1}(t)\le ...\le \l_{N}(t)$. From the spectral
theory of Jacobi operators \cite{vM} we have that if all $a_k(t)\ne
0$ for some $t$, then $\l_{-N}(t)< \l_{-N+1}(t)< ...< \l_{N}(t)$ and
perturbation theory gives us that any of these functions is analytic
in some neighborhood of $t$. \BBox

%\no {\bf Acknowledgments.}
% \footnotesize
%Some parts of this paper were written at the Universit´e Bordeaux 1;
%A.K. would like to thank to the Laboratoire deM´ecanique Physique
%(LMP) of the Universit´e Bordeaux 1 for the hospitality.

\end{document}